\documentclass[11pt]{article}

\usepackage{amsmath}
\usepackage{amsfonts}
\usepackage{amssymb}
\usepackage{latexsym}

\usepackage{filemod}

\newtheorem{theo}{Theorem}[section]
\newtheorem{lemma}[theo]{Lemma}
\newtheorem{coro}[theo]{Corollary}

\def\N{\mathbb N}
\def\R{\mathbb R}

\def\proof{{\bf Proof }}
\def\proofbox{\mbox{\bf Q.E.D.}\\}
\newcommand{\acos}{\mathrm{arccos\,}}
\begin{document}

\author{K\'aroly B\"or\"oczky\footnote{Roland E\"otv\"os University, P\'azm\'any P\'eter s\'et\'any 1/C, H-1117 Budapest, Hungary, boroczky@cs.elte.hu},
K\'aroly J. B\"or\"oczky\footnote{Alfr\'ed R\'enyi Institute of Mathematics, Hungarian Academy
  of Sciences, Re\'altanoda u. 13-15., H-1053 Budapest, Hungary, carlos@renyi.hu and Central European University, N\'ador u. 9., H-1051, Budapest, Hungary, Research is supported in parts by NKFIH grants 109789, 121649 and 116451.},
Alexey Glazyrin \footnote{The University of Texas Rio Grande Valley, School of Mathematical \& Statistical Sciences, One West University Blvd, Brownsville, Texas, USA, Alexey.Glazyrin@utrgv.edu, Research is partially supported by NSF grant DMS-1400876},
\'Agnes Kov\'acs \footnote{Roland E\"otv\"os University, P\'azm\'any P\'eter s\'et\'any 1/C, H-1117 Budapest, Hungary, kovacsa@cs.elte.hu}}

\title{Stability of the simplex bound for packings by equal spherical caps determined by simplicial regular polytopes}

\maketitle

MSC2010 Subject class: 52C17

\begin{abstract}
It is well known that the vertices of any simplicial regular polytope in $\R^d$ determine an optimal packing
of equal spherical balls in $S^{d-1}$. We prove a stability version of optimal order of this result.
\end{abstract}

\section{Introduction}

Euclidean regular polytopes are in the center of scientific studies since the Antiquity (see
P. McMullen, E. Schulte \cite{McS02} or H.S.M. Coxeter \cite{Cox73}). Packings of equal balls in spaces of constant curvature have been
investigated rather intensively since the
middle of the 20th century (see K. Bezdek \cite{Bez10} and G. Fejes T\'oth \cite{FTG04}).
In this paper, we focus on  packings of equal spherical balls (see
J.H.~Conway, N.J.A.~Sloane \cite{CoS98}, T. Ericson, V. Zinoviev \cite{ErZ01} and O. Musin \cite{Mus08}) that are related to some Euclidean simplicial regular polytope $P$ with its $f_0(P)$ vertices being on $S^{d-1}$, $d\geq 3$. We write $\varphi_P$ to denote the acute angle satisfying that
edge length of $P$ is $2\sin \varphi_P$. We note that the simplicial regular polytopes in $\R^d$, $d\geq 3$, are the regular simplex and crosspolytope in all dimensions, and in addition the icosahedron in $\R^3$ and the $600$-cell in $\R^4$ (the latter has Schl\"afli symbol $(3,3,5)$). The corresponding data is summarized in the following table.

$$
\begin{array}{|l|c|c|}
\hline
&&\\
\mbox{Regular Polytope }P&f_0(P)&\varphi_P\\
\hline \hline
\mbox{simplex in $\R^d$} & d+1 & \frac12\,\acos\frac{-1}{d}\\[0.5ex]
\hline
\mbox{crosspolytope in $\R^d$} & 2d & \frac\pi{4}\\[0.5ex]
\hline
\mbox{icosahedron in $\R^3$} & 12 & \frac12\,\acos\frac{1}{\sqrt{5}}\\[0.5ex]
\hline
\mbox{$600$-cell in $\R^4$} & 120 & \frac\pi{10}\\[0.5ex]
% && \\
\hline
\end{array}
$$

\noindent{\bf Theorem A }{\it If $P$ is a simplicial regular polytope in $\R^d$ having its vertices on $S^{d-1}$, $d\geq 3$, then the  vertices are centers of an optimal packing
of  equal spherical balls of radius $\varphi_P$ on $S^{d-1}$.}\\

Theorem A is due to Jung \cite{Jun01} if $P$ is a regular simplex. For the case of a regular crosspolytope, the statement of Theorem A was proposed as a problem by H. Davenport and Gy. Haj\'os \cite{DaH51}.
Numerous solutions arrived in a relatively short time;
namely, the ones by J. Acz\'el \cite{Acz52} and by T. Szele \cite{Sze52}
and the unpublished ones due to M. Bogn\'ar,
\'A. Cs\'asz\'ar, T. K\H ov\'ari and I. Vincze.
Independently, R.A. Rankin \cite{Ran55} solved the case of crosspolytopes. There are two more simplical regular polytopes. The case of icosahedron was handled by L. Fejes T\'oth \cite{FTL49} (see, say, \cite{FTL64} or \cite{FTL72}), and the case of the $600$-cell is due to K. B\"or\"oczky \cite{Bor78}.
All these arguments yield (explicitly or hidden) also the uniqueness of the optimal configuration up to orthogonal transformations.
For the case of the 600-cell, N.N. Andreev \cite{And99} provided an argument for optimality based on the linear programming bound in coding theory. The proof of uniqueness via the linear programming bound was given by P. Boyvalenkov and D. Danev \cite{BD01}.

In this paper, we provide a stability version of Theorem~A of optimal order. For $u,v\in S^{d-1}$, we write $\delta(u,v)\in[0,\pi]$ to denote the spherical (geodesic) distance
of $u$ and $v$, which is just their angle as vectors in $\R^d$.

\begin{theo}
\label{simpl-pol-stability}
Let $P$ be a simplicial regular polytope in $\R^d$ having its vertices on $S^{d-1}$, $d\geq 3$. For suitable $\varepsilon_P,c_P>0$, if $x_1,\ldots,x_{k}\in S^{d-1}$ are centers of non-overlapping spherical balls of radius at least $\varphi_P-\varepsilon$ for $\varepsilon\in[0,\varepsilon_P)$ and $k\geq f_0(P)$, then $k=f_0(P)$,
and  there exists a $\Phi\in O(d)$, such that  for any $x_i$ one finds a vertex $v$ of $P$ satisfying $\delta(x_i,\Phi v)\leq c_P\varepsilon$.
\end{theo}

We even provide explicit expressions for $\varepsilon_P$ and $c_P$. If $P$ is a $d$-simplex or a $d$-crosspolytope, then $c_P$ is of polynomial growth in $d$ ($c_P=9d^{3.5}$ if $P$ is a $d$-simplex, and $c_P=96d^3$ if $P$ is a $d$-crosspolytope).

Concerning notation, if  $p\in S^{d-1}$ and $\varphi\in(0,\pi/2)$, then we write $B(p,\varphi)$ for the spherical ball of center $p$ and radius $\varphi$. 
When working in $\R^d$, we write either $|X|$ or
 ${\cal H}^{d-1}(X)$ to denote the $(d-1)$-dimensional Hausdorff-measure of $X$. For $x_1,\ldots,x_k\in \R^d$, their convex hull, linear hull and affine hull in $\R^d$ are denoted by $[x_1,\ldots,x_k]$, ${\rm lin}\{x_1,\ldots,x_k\}$ and
${\rm aff}\{x_1,\ldots,x_k\}$, respectively.
For $x,y\in\R^d$, we write $\langle x, y\rangle$ to denote the scalar product, and $\|x\|$ to denote the Euclidean norm.
As usual, ${\rm int}\,K$ stands for the interior of $K\subset \R^d$.

The paper uses various tools to establish Theorem~\ref{simpl-pol-stability}. Only elementary linear algebra is needed for the case of a simplex, the linear programming bound is used for the case of a crosspolytope, and the simplex bound is applied to the icosahedron and the 600-cell.

Concerning the structure of the paper, Section~\ref{secsimplex} and Section~\ref{seccrosslin} handle the cases of the simplex and
 the crosspolytope, respectively, and Section~\ref{seclinear} in between reviews the linear programming bound used for the case of crosspolytopes. Results in these sections will be used also to handle the cases of the icosahedron in Section~\ref{seciso} and
the $600$-cell in Section~\ref{sec600}, as well. After reviewing the Delone and Dirichlet-Voronoi cell decompositions and the corresponding simplex bound in
Section~\ref{seccells}, and verifying some volume estimates in Section~\ref{secvol}, Theorem~\ref{simpl-pol-stability} is proved in
Section~\ref{seciso} and Section~\ref{sec600}  in the cases of the icosahedron and the $600$-cell, respectively.

\section{Some simple preparatory statements }
\label{secprepare}

The following statement will play a key role in the arguments for the cases of simplices and crosspolytopes
of Theorem~\ref{simpl-pol-stability}.

\begin{lemma}
\label{almostorthogonal}
Let $n\geq 2$ and $0\leq \eta<\frac1{n-1}$. If $u_1,\ldots,u_n\in S^{n-1}$ satisfy that
$|\langle u_i, u_j\rangle|\leq \eta$ for $i\neq j$, then there exists an orthonormal
basis $v_1,\ldots,v_n$ of $\R^n$ such that
${\rm lin}\{u_i,\ldots,u_n\}={\rm lin}\{v_i,\ldots,v_n\}$ and
$\langle u_i, v_i\rangle>0$
for $i=1,\ldots,n$, and
\begin{equation}
\label{uivjeta}
|\langle u_i, v_j\rangle|\leq \frac{\eta}{1-(n-2)\eta}\mbox{ \ for $i\neq j$.}
\end{equation}
Moreover, $\delta(u_i,v_i)\leq 2n\eta$ holds for $i=1,\ldots,n$ provided that $\eta<\frac1{2n}$.
\end{lemma}
\proof We prove the lemma by induction on $n$ where the case $n=2$ readily holds. Therefore, we assume that $n\geq 3$, and the lemma holds in $\R^{n-1}$.

Let $v_n=u_n$.
For $i=1,\ldots,n-1$, let $u_i=w_i+t_iv_n$ for $w_i\in v_n^\bot$ and $t_i\in\R$. It follows that
$|t_i|\leq \eta$ and $\|w_i\|=(1-t_i^2)^{\frac{1}2}\geq  (1-\eta^2)^{\frac{1}2}$
for $i=1,\ldots,n-1$, and we define $\bar{w}_i=w_i/\|w_i\|\in S^{n-1}$. We observe that
if $1\leq i<j\leq n-1$, then
$$
|\langle \bar{w}_i, \bar{w}_j\rangle|=
\frac{|\langle w_i, w_j\rangle|}{(1-t_i^2)^{\frac{1}2}(1-t_j^2)^{\frac{1}2}}\leq
\frac{|\langle u_i, u_j\rangle|+|t_it_j|}{1-\eta^2}\leq \frac{\eta+\eta^2}{1-\eta^2}
=\frac{\eta}{1-\eta}.
$$
As $\bar{\eta}=\frac{\eta}{1-\eta}<\frac1{n-2}$ follows from $\eta<\frac1{n-1}$, we may apply the induction hypothesis to $\bar{w}_1,\ldots,\bar{w}_{n-1}$ and $\bar{\eta}$. We obtain an orthonormal basis
$v_1,\ldots,v_{d-1}$ for $v_n^\bot$ such that
${\rm lin}\{\bar{w}_i,\ldots,\bar{w}_{n-1}\}={\rm lin}\{v_i,\ldots,v_{n-1}\}$ and
$\langle \bar{w}_i,v_i\rangle > 0$
for $i=1,\ldots,n-1$, and
$$
|\langle \bar{w}_i, v_j\rangle|\leq \frac{\bar{\eta}}{1-(n-3)\bar{\eta}}=\frac{\eta}{1-(n-2)\eta}\mbox{ \ for $i\neq j$.}
$$

If $1\leq i\leq n-1$ then $\langle u_n, v_i\rangle=\langle v_n, v_i\rangle=0$ and
$|\langle u_i, v_n\rangle|=|t_i|\leq \eta$. However if $i\neq j$ for $i,j\in\{1,\ldots,n-1\}$, then
$$
|\langle u_i, v_j\rangle|=|\langle (1-t_i^2)^{\frac{1}2}\bar{w}_i+t_iv_n, v_j\rangle|\leq
|\langle \bar{w}_i, v_j\rangle|\leq \frac{\eta}{1-(n-2)\eta}.
$$
Therefore, we have verified (\ref{uivjeta}), and we readily have
${\rm lin}\{u_i,\ldots,u_n\}={\rm lin}\{v_i,\ldots,v_n\}$ for $i=1,\ldots,n$ by construction.

Finally, for the estimate $\delta(u_i,v_i)$ if $\eta<\frac1{2n}$ and $i=1,\ldots,n$, we observe that
$|\langle u_i, v_j\rangle|<2\eta$ provided $j\neq i$.
It follows from $\|u_i\|=1$ and $\langle u_i,v_i\rangle>0$ that
$$
0\leq\langle v_i-u_i,v_i\rangle=1-\sqrt{1-\sum_{j\neq i}\langle u_i,v_j\rangle^2}\leq
\sum_{j\neq i}\langle u_i,v_j\rangle^2
\leq (n-1)4\eta^2<2\eta.
$$
In particular,
$$
\|v_i-u_i\|=\sqrt{\sum_{j=1}^n\langle v_i-u_i,v_j\rangle^2}<
\sqrt{n4\eta^2}=2\sqrt{n}\,\eta,
$$
and hence $\delta(u_i,v_i)<2n\eta$.
\proofbox

The following Lemma~\ref{inhemisphere} and its consequence Corollary~\ref{d+2pack} are due to
R.A. Rankin \cite{Ran55}, and will be used, say, for simplices.

\begin{lemma}
\label{inhemisphere}
If $u_1,\ldots,u_{d+1}\in S^{d-1}$, $d\geq 2$, are contained in a closed hemisphere, then there exist $i$ and $j$,
$1\leq i<j\leq d+1$, such that $\langle u_i, u_j\rangle\geq 0$.
\end{lemma}
\proof We prove the statement by induction on $d$ where the case $d=2$ readily holds. If $d\geq 3$, then
we may assume that $\langle u_i, u_j\rangle\leq 0$ if $1\leq i<j\leq d+1$.
Let
$v\in S^{n-1}$ such that $\langle v, u_i\rangle\geq 0$ for $i=1,\ldots,d+1$, and hence $u_i=w_i+\lambda_iv$
for $i=1,\ldots,d+1$ where $w_i\in v^\bot$ and $\lambda_i\geq 0$. If $u_i=v$ for some
$i\in\{1,\ldots,d+1\}$, then $\langle u_j, u_i \rangle=0$ for $j\neq i$, thus we are done.
Otherwise $w_i\neq o$ for $i=1,\ldots,d+1$. If $i=1,\ldots,d$, then
$$
0\geq \langle u_{d+1}, u_i\rangle =\langle w_{d+1}, w_i\rangle +\lambda_{d+1}\cdot \lambda_i\geq \langle w_{d+1}, w_i\rangle,
$$
therefore, the induction hypothesis applied to
$\frac{w_1}{\|w_1\|},\ldots,\frac{w_d}{\|w_d\|}\in v^\bot\cap S^{d-1}$ yields
$\langle w_i, w_j\rangle \geq 0$ for some $1\leq i<j\leq d$, and hence
$\langle u_i, u_j\rangle\geq 0$.
\proofbox

\begin{coro}
\label{d+2pack}
If  $k\geq d+2$, $d\geq 2$, and $u_1,\ldots,u_k\in S^{d-1}$, then there exist $i$ and $j$,
$1\leq i<j\leq d+1$, such that $\langle u_i, u_j\rangle\geq 0$.
\end{coro}

\section{The proof of Theorem~\ref{simpl-pol-stability} in the case of Simplices }
\label{secsimplex}

Theorem~\ref{simplex-unit-vectors} covers the case of regular simplex of Theorem~\ref{simpl-pol-stability}.

\begin{theo}
\label{simplex-unit-vectors}
 If $u_0,\ldots,u_d\in S^{d-1}$ satisfy
$\delta(u_i,u_j)\geq \acos\frac{-1}d-2\varepsilon$ for $\varepsilon\in[0,\varepsilon_d)$ and $0\leq i<j\leq d$, $d\geq 2$,
 then there exists a regular simplex $[v_0,\ldots,v_d]$ with $v_0,\ldots,v_d\in S^{d-1}$
 such that
$\delta(u_i,v_i)\leq c_d\varepsilon$ for $i=0,\ldots, d$ where $c_d=9d^{3.5}$ and $\varepsilon_d=1/c_d$.
\end{theo}
{\bf Remark } If $d=2$, then one may even choose $c_2=3$ and  $\varepsilon_2=\frac{\pi}{12}$.\\
\proof We first handle the case $d=2$, because this case is much more elementary. We define
$\varepsilon_2$ to be $\frac{\pi}{12}=\frac12(\frac{2\pi}3-\frac{\pi}2)$. Thus
$\acos\frac{-1}2=\frac{2\pi}3$ and $\varepsilon<\varepsilon_2$ yield that no closed semicircle contains $u_0,u_1,u_2$, and hence the sum of the three angles
of type $\delta(u_i,u_j)$ is $2\pi$.
We may assume that
$\delta(u_0,u_1)\leq \delta(u_0,u_2)\leq \delta(u_1,u_2)$, and hence
\begin{equation}
\label{2dimeps}
\frac{2\pi}3-2\varepsilon\leq \delta(u_0,u_1)\leq \frac{2\pi}3
\leq \delta(u_1,u_2)\leq \frac{2\pi}3+4\varepsilon.
\end{equation}
We choose $v_1,v_2,v_3\in S^1$ that are vertices of a regular triangle, and
$$
\delta(u_0,v_0)=\delta(u_1,v_1)\leq \varepsilon.
$$
We deduce from (\ref{2dimeps}) that $\delta(u_2,v_2)\leq 3\varepsilon$, thus we may choose $c_2$ to be $3$.

Turning to the case $d\geq 3$, let
$$
0<\varepsilon<\frac1{9d^{3.5}}.
$$
 If $0\leq i<j\leq d$, then we have
\begin{eqnarray}
\nonumber
\|u_i-u_j\|^2&=&2-2\cos\delta(u_i,u_j)\geq 2+2\left(\frac{\cos 2\varepsilon}d -
\frac{\sqrt{d^2-1}}d\cdot \sin 2\varepsilon\right)\\
\label{uiujlow}
&>&2+2\left(\frac{1-2\varepsilon}d -
 2\varepsilon\right)
>\frac{2(d+1)}d-6\varepsilon.
\end{eqnarray}
Using (\ref{uiujlow}) and the estimate
$$
(d+1)^2=\left\|\sum_{i=0}^du_i\right\|^2
+\sum_{0\leq i<j\leq d}\|u_i-u_j\|^2\geq \sum_{0\leq i<j\leq d}\|u_i-u_j\|^2,
$$
we deduce for any $i<j$ the upper bound
$$
\|u_i-u_j\|^2<\frac{2(d+1)}d+3d(d+1)\varepsilon.
$$
In particular, if  $i< j$, then
\begin{equation}
\label{xycondsimplex}
\frac{-1}d-\frac32\,d(d+1)\varepsilon\leq\langle u_i,u_j\rangle\leq \frac{-1}d+ 2\varepsilon.
\end{equation}

We  embed $\R^d$ into $\R^{d+1}$ as $\R^d=e^\bot$ for suitable
$e\in S^d\subset \R^{d+1}$. For $i=0,\ldots,d$, we define
$$
w_i=\sqrt{\frac1{d+1}}\,e+\sqrt{\frac{d}{d+1}}\,u_i\in S^d,
$$
and hence (\ref{xycondsimplex}) yields that if $i\neq j$, then
$$
|\langle w_i,w_j\rangle|=\left|\frac1{d+1}+\frac{d}{d+1}\langle u_i,u_j\rangle\right|=
\frac{d}{d+1}\left|\frac1d+\langle u_i,u_j\rangle\right|\leq \frac32\,d^2\varepsilon.
$$
Since $\frac32\,d^2\varepsilon<\frac1{2(d+1)}$, Lemma~\ref{almostorthogonal} can be applied, and hence there exists an orthonormal basis $q_0,\ldots,q_d$
of $\R^{d+1}$ such that $\delta(w_i,q_i)\leq 3(d+1)d^2\varepsilon$ holds for $i=0,\ldots,d$.
We define $q=\sum_{i=0}^d\frac1{\sqrt{d+1}}\,q_i$ and deduce that $q\in S^d$.

Since for any $i=0,\ldots,d$, we have $\langle e,w_i\rangle=\frac1{\sqrt{d+1}}$ and $\delta(w_i,q_i)\leq 3(d+1)d^2\varepsilon$, it follows from
$|\cos(\alpha+\beta)-\cos\alpha|\leq|\beta|$  for $\alpha,\beta\in\R$ that
$\left|\langle e,q_i\rangle-\frac1{\sqrt{d+1}}\right|\leq 3(d+1)d^2\varepsilon$,
and hence $|\langle e-q,q_i\rangle|\leq 3(d+1)d^2\varepsilon$.
We deduce that
$$
\|e-q\|\leq 3(d+1)^{\frac32}d^2\varepsilon.
$$
Let $A\in{\rm O}(d+1)$ be the identity if $e=q$, and be the rotation around the linear $(d-1)$-space of $\R^{d+1}$ orthogonal to ${\rm lin}\{e,q\}$ with $Aq=e$ if $e\neq q$. It follows that $\|Au-u\|\leq \|e-q\|$ for $u\in S^d$. For each $i=0,\ldots,q$, $\bar{q}_i=Aq_i\in S^d$ satisfies
$\|\bar{q}_i-q_i\|\leq\|e-q\|\leq 3(d+1)^{\frac32}d^2\varepsilon$ and combining the last estimate
with $\delta(w_i,q_i)\leq 3(d+1)d^2\varepsilon\leq\frac 32 (d+1)^{\frac 32}d^2\varepsilon$   yields
\begin{equation}
\label{eqdist}
\|w_i-\bar{q}_i\|\leq \frac 92 (d+1)^{\frac32}d^2\varepsilon.
\end{equation}
As $Aq=e$, we also have that
$\langle \bar{q}_i,e\rangle=\sqrt{\frac1{d+1}}=\langle w_i,e\rangle$ for $i=0,\ldots,q$.
Therefore,
$$
v_i=\sqrt{\frac{d+1}{d}}\left(\bar{q}_i-\sqrt{\frac1{d+1}}\,e\right)\in e^\bot\cap S^d=S^{d-1}
$$
for $i=0,\ldots,q$, $[v_0,\ldots,v_d]$ is a regular $d$-simplex, and
$$
\|v_i-u_i\|=\sqrt{\frac{d+1}{d}}\cdot \|\bar{q}_i-w_i\|\leq \frac 92 (d+1)^2d^{\frac32}\varepsilon\leq 8 d^{3.5}\varepsilon\leq \frac 89
$$
for $i=0,\ldots,q$ where we used $d\geq 3$ at the last estimate. Using that $2 \arcsin \frac t 2\leq \frac 98 t$ for any $t\in[0,\frac 89]$, we conclude that $\delta(v_i,u_i)=2\arcsin \frac {\|v_i-u_i\|} 2 \leq \frac 98 \|v_i-u_i\|\leq 9d^{3.5}\varepsilon$ for $i=0,\ldots,q$. \proofbox

\section{The linear programming bound}
\label{seclinear}

Let $d\geq 2$. The presentation about the linear programming bound for sphere packings on $S^{d-1}$ in this section is based on T. Ericson, V. Zinoviev \cite[Chapter~2]{ErZ01}.
A central role in the theory is played by certain real Gegenbauer polynomials $Q_i$, $i\in\N$, in one variable where each $Q_i$ is of degree $i$, and satisfies the following recursion:
\begin{eqnarray*}
Q_0(t)&=&1\\
Q_1(t)&=&t\\
Q_2(t)&=&\frac{dt^2-1}{d-1}\\
(i+d-2)Q_{i+1}(t)&=&(2i+d-2)tQ_i(t)-iQ_{i-1}(t) \mbox{ \ for $i\geq 2$}.
\end{eqnarray*}
We do not signal the dependence of $Q_i$ on $d$ because the original notation for the Gegenbaur polynomial is
$Q_i=Q_i^{(\alpha)}$ for $\alpha=\frac{d-2}2$ as
$$
\int_{-1}^1Q_i(t)Q_j(t)(1-t^2)^{\frac{d-3}2}\,dt=0\mbox{ \ \ if $i\neq j$}.
$$
Actually, $Q_i$ is normalized in a way such that $Q_i(1)=1$ for $i\in \N$.

The basis of our considerations is the following version of the linear programming bound, which is contained in the proof of Theorem~2.3.1 in \cite{ErZ01}. We write $|X|$ to denote the cardinality of a finite set $X$.

\begin{theo}
\label{linprogbound}
For $d\geq 2$, if $f=f_0Q_0+f_1Q_1+\ldots+f_kQ_k$ for $k\geq 1$, $f_0>0$ and
 $f_1,\ldots,f_k\geq 0$, then any finite $X\subset S^{d-1}$ satisfies
\begin{equation}
\label{linprogbound0}
|X|f(1)+\sum_{x,y\in X\atop x\neq y}f(\langle x,y\rangle)\geq |X|^2 f_0.
\end{equation}
\end{theo}
{\bf Remark} The classical linear programming bound is a consequence; namely, if in addition, $f(t)\leq 0$ for
fixed $s\in(-1,1)$ and variable $t\in[-1,s]$, then
\begin{equation}
\label{linprogbound1}
|X|\leq f(1)/ f_0.
\end{equation}

If we have equality in (\ref{linprogbound1}), then (\ref{linprogbound0}) shows that all values
$\langle x,y\rangle$ for $x\neq y$, $x,y\in X$ are roots of $f$.

As an example, let $X\subset S^{d-1}$ be the centers for a packing of spherical balls of radius $\frac{\pi}4$,
and hence
$\langle x,y\rangle\leq 0$ for $x,y\in X$ with $x\neq y$. The polynomial
$$
f(t)=t(t+1)=f_0Q_0+f_1Q_1+f_2Q_2
$$
satisfies $f(t)\leq 0$ for $t\in[-1,0]$ and
$$
f_0=\frac1d,\mbox{ \ }f_1=1,\mbox{ \ }f_2=1-\frac1d,\mbox{ \ }f(1)=2,
$$
therefore, (\ref{linprogbound1}) yields $|X|\leq 2d$.

Next we quantify the obvious statement that for any packing of $m$  spherical balls of radius $r$ on $S^{n-1}$, if $r$ is close to $\frac{\pi}4$ then $m\leq 2n$.

\begin{lemma}
\label{almostcross}
If $Y\subset S^{n-1}$, $n\geq 2$, satisfies that
$\langle x,y\rangle<\frac1{2n^2-n}$ for $x,y\in Y$ with $x\neq y$, then $|Y|\leq 2n$.
\end{lemma}
\proof Let $s=\max\{\langle x,y\rangle:\,x,y\in Y\mbox{ and }x\neq y\}<\frac1{2n^2-n}$. We consider the polynomial
$$
f(t)=(t+1)(t-s)=f_0Q_0+f_1Q_1+f_2Q_2
$$
where $f(t)\leq 0$ for $t\in[-1,s]$ and
$$
f_0=\frac1n-s,\mbox{ \ }f_1=1-s,\mbox{ \ }f_2=1-\frac1n,\mbox{ \ }f(1)=2(1-s).
$$
We deduce from the linear programming bound (\ref{linprogbound1}) and $s<\frac1{2n^2-n}$ that
$$
|Y|\leq \frac{2n(1-s)}{1-ns}=2n+\frac{2n(n-1)s}{1-ns}<2n+1. \mbox{ \ }\proofbox
$$

The linear programming bound could have been used in the case of simplex to prove (\ref{xycondsimplex}). However, this could be proved easily by elementary arguments, as well.

The linear programming bound can be also used to prove the optimality of the icosahedron and the 600-cell however the corresponding polynomials are more complicated. Say, in the case of 600-cell, the polynomial is of degree $17$ and $f_{12}=f_{13}=0$ according to N.N. Andreev \cite{And99}. Therefore we use volume estimates to handle the cases of
the icosahedron and the 600-cell.

\section{The proof of Theorem~\ref{simpl-pol-stability} in the case of Crosspolytopes }
\label{seccrosslin}

Let $X\subset S^{d-1}$ be the centers for a packing of at least $2d$ spherical balls of radius $\frac{\pi}4-\varepsilon$,
$0<\varepsilon <\frac1{64d^4}$, and hence
$\langle x,y\rangle\leq s$ for $x,y\in X$ with $x\neq y$ and
$$
s=\sin 2\varepsilon < 2\varepsilon <\frac1{32d^4}.
$$
 We deduce from Lemma~\ref{almostcross} that
$$
|X|=2d.
$$

We consider the polynomial
$$
f(t)=(t+1)(t-s)=f_0Q_0+f_1Q_1+f_2Q_2
$$
where $f(t)\leq 0$ for $t\in[-1,s]$ and
$$
f_0=\frac1d-s,\mbox{ \ }f_1=1-s,\mbox{ \ }f_2=1-\frac1d,\mbox{ \ }f(1)=2(1-s).
$$
It follows from (\ref{linprogbound0}) and $f(t)\leq 0$ for $t\in[-1,s]$ that
if $x,y\in X$ with $x\neq y$, then
\begin{equation}
\label{fxycross}
f(\langle x,y\rangle)\geq |X|^2 f_0-|X|f(1)=4d^2\left(\frac1d-s\right)-4d(1-s)=-4d(d-1)s.
\end{equation}

Since $t-s\leq \frac{-1}2$ if $t\leq \frac{-1}2$ and $t+1\geq \frac{1}2$ if $t\geq \frac{-1}2$, we have
$$
f(t)\leq -\frac12\min\left\{|t+1|,|t-s|\right\}\mbox{ \ \ for $t\in[-1,s]$}.
$$
We deduce from (\ref{fxycross}) that if $x,y\in X$ with $x\neq y$, then
$$
\min\{\langle x,y\rangle+1, s-\langle x,y\rangle\}\leq 8d(d-1)s,
$$
or in other words,
\begin{equation}
\label{xycross}
\begin{array}{rrcl}
\mbox{either }&-1\leq&\langle x,y\rangle&\leq -1+\frac 1 {4d^2}<\frac{-3}4\\[1ex]
\mbox{or }&-8d(d-1)s\leq&\langle x,y\rangle&\leq s<\frac1{32}.
\end{array}
\end{equation}
We define
\begin{equation}
\label{etadef}
\eta=8d(d-1)s<\frac1{4d^2}.
\end{equation}

We claim that for every $x\in X$
\begin{equation}
\label{xoppositecross}
\mbox{there exists a unique $y\in X$ such that }\langle x,y\rangle\leq \frac{-3}4,
\end{equation}
which we call the element of $X$ opposite to $x$. For any $y\in X$, we write
$\bar{y}$ to denote its projection into $x^\bot$, and if $y\neq\pm x$, then
we set $y^*=\bar{y}/\|\bar{y}\|$.

The first step towards (\ref{xoppositecross}) is to show that if $y,z\in X$, then
\begin{equation}
\label{xyzcross}
\langle x,y\rangle\leq \frac{-3}4\mbox{ and }\langle x,z\rangle\leq \frac{-3}4
\mbox{ \ yield } y=z.
\end{equation}
Since
$\|\bar{y}\|=\sqrt{1-\langle x,y\rangle^2}<\sqrt{\frac12}$ and similarly
$\|\bar{z}\|<\sqrt{\frac12}$, we have
$$
\langle y,z\rangle=\langle x,y\rangle\langle x,z\rangle+\langle \bar{y},\bar{z}\rangle>
\frac{9}{16}-\frac12=\frac1{16},
$$
which proves  $\langle y,z\rangle=1$ by (\ref{xycross}), and in turn verifies (\ref{xyzcross}).

Next, set $\widetilde{X}=\{y\in X:\,|\langle x,y\rangle|\leq \eta\}$. For (\ref{xoppositecross}), it is sufficient to verify that
\begin{equation}
\label{Xtildecard}
|\widetilde{X}|\leq 2(d-1).
\end{equation}
For $y_1,y_2\in \widetilde{X}$, we have $y_i=\bar{y}_i+p_ix$
for $i=1,2$ where $p_i\in[-\eta,\eta]$.
In particular, $\|\bar{y}_i\|=(1-p_i^2)^{\frac{1}2}\geq  (1-\eta^2)^{\frac{1}2}$, and hence
$$
\langle y_1^*,y_2^*\rangle=
\frac{\langle \bar{y}_1,\bar{y}_2\rangle}{(1-p_1^2)^{\frac{1}2}(1-p_2^2)^{\frac{1}2}}=
\frac{\langle y_1, y_2\rangle-p_1p_2}{(1-p_1^2)^{\frac{1}2}(1-p_2^2)^{\frac{1}2}}
\leq \frac{\eta+\eta^2}{1-\eta^2}
=\frac{\eta}{1-\eta}<2\eta.
$$
Since $2\eta<\frac1{2d^2}$, Lemma~\ref{almostcross} with $n=d-1$ yields
(\ref{Xtildecard}), and in turn (\ref{xoppositecross}).

We deduce from (\ref{xoppositecross}) that $X$ can be divided into $d$ pairs of opposite vectors. Choosing one unit vector from each pair, we obtain $x_1,\ldots,x_d\in X$ such that $|\langle x_i,x_j\rangle|\leq \eta$ for $i\neq j$. It follows from Lemma~\ref{almostorthogonal} that for every such $d$-tuple $x_1,\ldots,x_d\in X$ there exists an orthonormal
basis $v_1,\ldots,v_d$ of $\R^d$ such that
${\rm lin}\{x_i,\ldots,x_d\}={\rm lin}\{v_i,\ldots,v_d\}$ and
$\delta(x_i,v_i)\leq 2d\eta$
for $i=1,\ldots,d$.

We claim that if $x,y\in X$ are opposite vectors, then
\begin{equation}
\label{oppositeantipodal}
\delta(y,-x)\leq 4d\eta.
\end{equation}
We choose $x_2,\ldots,x_d\in X$ representatives from the other $d-1$ opposite pairs, and let $v$ be the unit vector
orthogonal to ${\rm lin}\{x_2,\ldots,x_d\}$ with $\langle x,v\rangle>0$.
Taking $x=x_1$ and considering the approximating orthonormal basis $v_1,\ldots,v_d$ for
this $x_1,\ldots,x_d$, we deduce that $v=v_1$, and hence $\delta(x,v)\leq 2d\eta$.
Similarly, taking $y=x_1$, we have $v_1=-v$ for the approximating orthonormal basis, thus
$\delta(y,-v)\leq 2d\eta$. In turn, we conclude (\ref{oppositeantipodal}) by the triangle inequality.

Finally, we fix representatives $u_1,\ldots,u_d$ from each of the $d$ pairs of opposite vectors, and hence
there exists an orthonormal
basis $w_1,\ldots,w_d$ of $\R^d$ such that
$\delta(u_i,w_i)\leq 2d\eta$
for $i=1,\ldots,d$. We write $u_{i+d}$ to denote the vector of $X$ opposite to $u_i$, $i=1,\ldots,d$, and hence
$\delta(u_{i+d},-u_i)\leq 4d\eta$ according to (\ref{oppositeantipodal}). Therefore,
$$
\delta(u_{i+d},-w_i)\leq \delta(u_{i+d},-u_i)+\delta(-u_i,-w_i)\leq 4d\eta+2d\eta=6d\eta\leq
48d^3s\leq 96d^3\varepsilon.
$$
Therefore, $c_d=96d^3$ can be chosen for Theorem~\ref{simpl-pol-stability} in the case of crosspolytopes.

\section{Spherical Dirichlet-Voronoi and Delone cell decomposition}
\label{seccells}

For $v\in S^{d-1}$ and acute angle $\theta$, we write $B(v,\theta)$ to denote the spherical ball of center $v$ and radius $\theta$. For $u,v\in S^{d-1}$, $u\neq -v$, we write $\overline{uv}$ to denote the smaller geodesic arc connecting $u$ and $v$. We will frequently use the Spherical Law of Cosines: If $a,b,c$ are side lengths of a spherical triangle contained in an open hemisphere, and the opposite angles are $\alpha,\beta,\gamma$, respectively, then
\begin{equation}
\label{LawCosine}
\cos c=\cos a \cdot \cos b+\sin a\cdot \sin b\cdot \cos\gamma.
\end{equation}

A  set $C\subset \R^d$ is a convex cone if it is closed and $\alpha x+\beta y\in C$ for $\alpha,\beta\geq 0$ and
$x,y\in C$. If $C$ contains a half-line, then $M=C\cap S^{d-1}$ is called a spherically convex set whose dimension is one less than the Euclidean dimension of $C$. The relative interior of $M$ is the intersection of $S^{d-1}$ and the relative interior of $C$ with respect to ${\rm lin}\,C$. If the origin is a face of $C$ and $C$ is a polyhedron (namely, intersection of finitely many half-spaces) then $M$ is called a spherical polytope. In this case, the faces of $M$ are intersections of $S^{d-1}$ with the faces of $C$ different from the origin.

Let $x_1,\ldots,x_k\in S^{d-1}$ satisfy that each open hemisphere contains some of $x_1,\ldots,x_k$, and hence
$o\in{\rm int}\,P$ for $P=[x_1,\ldots,x_k]$. The radial projections of the facets of $P$ onto $S^{d-1}$ form the Delone
(or Delaunay) cell decomposition of $S^{d-1}$. We observe that if the distance of $o$ from ${\rm aff}\,F$ is $\varrho$ for a facet $F$, then $\arccos \varrho$ is the spherical radius of the spherical cap cut off by ${\rm aff}\,F$. We call
$\arccos \varrho$ the spherical circumradius of the corresponding Delone cell.

To define the other classical decomposition of $S^{d-1}$ corresponding to $x_1,\ldots,x_k$,
let
$$
D_i=\{u\in S^{d-1}:\,\delta(u,x_i)\leq \delta(u,x_j)\mbox{ for }j=1,\ldots,k\}
$$
for $i=1,\ldots,k$, which is the Dirichlet-Voronoi cell of $x_i$. The Dirichlet-Voronoi cells also form a cell decomposition of $S^{n-1}$ that is dual to the Delone cell decomposition by providing the following bijective correspondence between vertices of Dirichlet cells and Delone cells. If $v$ is a vertex of $D_i$, $i\in\{1,\ldots,k\}$, and $\delta(v,x_i)=\theta$, then $\delta(v,x_j)\geq\theta$ for all
$j=1,\ldots,k$, and points $x_j$ with $\delta(v,x_j)=\theta$ form the vertex set of a Delone  cell (see, say, K.J. B\"or\"oczky \cite{Bor04}). In addition, if $F$ is an $m$-dimensional face of some $D_i$, and
$p$ is the closest point of the $m$-dimensional great sphere $\Sigma$ of $F$, then there exists a $(d-1-m)$-dimensional face $G$ of the Delone cell complex contained in the $(d-1-m)$-dimensional great sphere $\Sigma'$ orthogonal to $\Sigma$ at $p$ whose vertices are all of distance $\delta(p,x_i)$ from $p$.

A simplex with ordered vertices $p_0,\ldots,p_{d-1}$ on $S^{d-1}$ is called an orthoscheme if
for $i=1,\ldots,d-2$, the $i$-dimensional great sphere through $p_0,\ldots,p_i$ is orthogonal to the
$(d-1-i)$-dimensional great sphere through $p_i,\ldots,p_{d-1}$.

For any face $F$ of a Dirichlet-Voronoi cell $D_i$, we write $q_i(F)$ to denote the point of $F$ closest to $x_i$.
It follows from the convexity of $F$ and the Spherical Law of Cosines that if $x\in F\backslash q_i(F)$, then
\begin{description}
\item{(a)} the angle
between the arcs $\overline{q_i(F),x_i}$ and $\overline{q_i(F),x}$ is at least $\frac{\pi}2$,
\item{(b)} and is actually exactly
$\frac{\pi}2$ if $q_i(F)$ lies in the relative interior of $F$.
\end{description}
For a Dirichlet-Voronoi cell $D_i$, we say that a sequence $(F_0,\ldots,F_{d-2})$ is a tower, if
 $F_j$ is a $j$-face of $D_i$, $j=0,\ldots,d-2$, and $F_j\subset F_l$ if $j<l$. In addition, $(F_0,\ldots,F_{d-2})$ is a
proper tower, if $q_i(F_j)\neq q_i(F_l)$ for $j<l$, and, in this case, we call the simplex $\Xi$ with ordered vertices
$x_i,q_i(F_{d-2}),\ldots,q_i(F_0),$ a quasi-orthoscheme. We observe that according to (b), a
 quasi-orthoscheme is an orthoscheme if each $q_i(F_j)$, $j=1,\ldots,d-2$, lies in the relative interior of $F_j$. Moreover, (a) yields that quasi-orthoschemes provide a triangulation of $S^{d-1}$ refining the Dirichlet-Voronoi cell decomposition.

For any $\varphi\in(0,\frac{\pi}2)$ and $i\geq 1$, we write $r_i(\varphi)\in(0,\frac{\pi}2)$ to denote the circumradius of the $i$-dimensional spherical
regular simplex of  edge length $2\varphi$. In particular, there exists a spherical triangle with equal sides $r_i(\varphi)$ enclosing the angle $\acos\frac{-1}i$ where the third side of the triangle is $2\varphi$. In addition, we define $r_\infty(\varphi)\in(0,\frac{\pi}2)$ in a way such that there exists a spherical triangle with equal sides
$r_\infty(\varphi)$ enclosing the right angle where the third side of the triangle is $2\varphi$. We have
$$
\varphi=r_1(\varphi)<\ldots<r_{d-1}(\varphi)<r_\infty(\varphi).
$$
It follows from (\ref{LawCosine}) that if $j=1,\ldots,d-1$, then
\begin{equation}
\label{rjdef}
\cos 2\varphi=\cos^2 r_j(\varepsilon)-\frac{\sin^2 r_j(\varepsilon)}{j}
\mbox{ \ and \ }\cos 2\varphi=\cos^2 r_\infty(\varepsilon),
\end{equation}
which in turn yields that
\begin{equation}
\label{rjdefsin}
\sin r_j(\varphi)=\sqrt{\frac{2j}{j+1}}\,\sin\varphi
\mbox{ \ and \ }\sin r_\infty(\varphi)=\sqrt{2}\,\sin\varphi.
\end{equation}

The following lemma is due to K. Boroczky \cite{Bor78}. We include the argument because the second statement is only implicit in \cite{Bor78}.

\begin{lemma}
\label{qiF}
Let $\varphi\in(0,\frac{\pi}2)$, and let $x_1,\ldots,x_k\in S^{d-1}$ satisfy that each open hemisphere contains some of
$x_1,\ldots,x_k$, and $\delta(x_i,x_j)\geq 2\varphi$ for $i\neq j$, and let $D_j$ be the Dirichlet-Voronoi cell of $x_j$. If $F$ is an $m$-dimensional face of certain $D_i$, then
\begin{description}
\item{(i)} $\delta(x_i,q_i(F))\geq r_{d-1-m}(\varphi)$;
\item{(ii)} and even $\delta(x_i,q_i(F))\geq r_{\infty}(\varphi)$ if $q_i(F)$ is not contained in the relative interior of $F$.
\end{description}
\end{lemma}
\proof Let $p$ be the closest to $x_i$ point of the $m$-dimensional great subsphere $\Sigma$ containing $F$, and let $I$ be the set of all indices $j$ such that $F$ is a face of $D_j$. In particular, all $x_j$ with $j\in I$ span the
$(d-1-m)$-dimensional great subsphere $\Sigma'$ passing through $p$ and perpendicular to $\Sigma$, and hence the cardinality of $I$ is at least $d-m$. It follows that for $\theta=\delta(x_i,p)\leq \delta(x_i,q_i(F))$, we have
$\theta=\delta(x_j,p)$ for $j\in I$. For $j\in I$, let $u_j$ be a unit vector tangent to the arc $\overline{p,x_j}$ at $p$,
and hence all $u_j$, $j\in I$, span the $(d-1-m)$-dimensional linear  subspace $L'$ tangent to $\Sigma'$ at $p$.
According to Jung's theorem (see also Lemma~\ref{simplex-unit-vectors}), there exist different $l,j\in I$ such that
$\delta(u_l,u_j)\leq \acos\frac{-1}{d-1-m}$. Since $\delta(x_l,p)=\delta(x_j,p)=\theta$, we deduce (i) from the Spherical Law of Cosines (\ref{LawCosine}).

Turning to (ii), we assume that $p$  is not contained in the relative interior of $F$. In this case, there exists an
$x_g\in S^{d-1}\backslash\Sigma'$ such that $0<\delta(x_g,p)\leq \theta$. Let $u_g\in S^{d-1}$ be a unit vector tangent to the arc $\overline{p,x_g}$ at $p$. We claim that there exist different $j,l\in I\cup\{g\}$ such that
\begin{equation}
\label{ujulI}
\langle u_j, u_l\rangle\geq 0.
\end{equation}
Let $L$ be the $m$-dimensional linear subspace $L$ tangent to $\Sigma$ at $p$, which is the orthogonal complement of $L'$ inside the tangent space to $S^{d-1}$ at $p$. Therefore,
there exist unit vectors $v\in L$ and $v'\in L'$ and a real number $t\in[0,\frac{\pi}2]$ such that
$u_g=v\cos t+v'\sin t$. If $\langle v', u_j\rangle< 0$ for all $j\in I$, then Lemma~\ref{inhemisphere} yields different $j,l\in I$ such that $\langle u_j, u_l\rangle\geq 0$. Otherwise there exists $j\in I$ such that $\langle v', u_j\rangle\geq 0$, and hence $\langle u_g, u_j\rangle\geq 0$, as well.

Using these $u_j$ and $u_l$ in (\ref{ujulI}), we apply the Spherical Law of Cosines (\ref{LawCosine}) to the triangle with vertices $p,x_j,x_l$ to obtain
$$
\cos 2\varphi\geq\cos\delta(x_j,x_l)\geq \cos\delta(p,x_j)\cdot \cos\delta(p,x_l)\geq \cos^2\theta.
$$
Therefore, $\theta\geq r_\infty(\varphi)$ by (\ref{rjdef}).
\proofbox

We fix a point $z_0\in S^{d-1}$, and for $0<t_1<\ldots<t_{d-1}<\frac{\pi}2$, we write $\Theta(t_1,\ldots,t_{d-1})$
to denote an orthoscheme with ordered vertices $z_0,z_1,\ldots,z_{d-1}$ such that
$\delta(z_0,z_i)=t_i$ for $i=1,\ldots,d-1$. We observe that the (spherical) diameter of
$\Theta(t_1,\ldots,t_{d-1})$ is $t_{d-1}$. For any $\varphi\in(0,t_1]$, we define
$$
\Delta(t_1,\ldots,t_{d-1})=\frac{|\Theta(t_1,\ldots,t_{d-1})\cap B(z_0,\varphi)|}
{|\Theta(t_1,\ldots,t_{d-1})|\cdot |B(z_0,\varphi)|},
$$
whose value does not depend on the choice of $\varphi\in(0,t_1]$. If $\Psi\subset z_0^\bot$
is the Euclidean convex polyhedral cone generated by the rays tangent to the arcs
$\overline{z_0,z_i}$ at $z_0$, $i=1,\ldots,d-1$, then
$$
\Delta(t_1,\ldots,t_{d-1})=\frac{{\cal H}^{d-2}(\Psi\cap S^{d-1})}
{|\Theta(t_1,\ldots,t_{d-1})|\cdot {\cal H}^{d-2}(S^{d-2})}.
$$
According to one of the core results
of K. Boroczky \cite{Bor78}, if $s_1<\ldots s_{d-1}<\frac{\pi}2$, and $t_i\leq s_i$ for $i=1,\ldots,d-1$, then
\begin{equation}
\label{orthoscheme-density-monotone}
\Delta(t_1,\ldots,t_{d-1})\geq \Delta(s_1,\ldots,s_{d-1}).
\end{equation}

We deduce from Lemma~\ref{qiF} and (\ref{orthoscheme-density-monotone}) the following estimate.

\begin{lemma}
\label{orthoscheme-density}
Let $\sigma\in(0,\frac{\pi}2)$, and let $x_1,\ldots,x_k\in S^{d-1}$, $d\geq 3$, satisfy that each open hemisphere contains some of
$x_1,\ldots,x_k$, and $\delta(x_i,x_j)\geq 2\sigma$ for $i\neq j$, and let $D_i$ be the Dirichlet-Voronoi cell of $x_i$. If
$\Xi$ is a quasi-orthoscheme associated to some $D_i$ and it is known that $\Xi$ is an orthoscheme, and the diameter of $\Xi$ is $R$, then
\begin{eqnarray}
\label{simplex-bound1}
\frac{|\Xi\cap B(x_i,\sigma)|}{|\Xi|\cdot |B(x_i,\sigma)|}
&\leq &\Delta(r_1(\sigma),\ldots,r_{d-2}(\sigma),R)\\
\label{simplex-bound2}
&\leq &\Delta(r_1(\sigma),\ldots,r_{d-2}(\sigma),r_{d-1}(\sigma)).
\end{eqnarray}
\end{lemma}

We note that the ideas in K. Boroczky \cite{Bor78} yield (\ref{simplex-bound2}) even if
the  quasi-orthoscheme $\Xi$ is not an orthoscheme, but they actually even imply the following stronger bound
in the low dimensions we are interested in.

\begin{lemma}
\label{quasi-orthoscheme-density}
Let $\sigma\in(0,\frac{\pi}2)$, and let $x_1,\ldots,x_k\in S^{d-1}$, $d=3,4$, satisfy that each open hemisphere contains some of
$x_1,\ldots,x_k$, and $\delta(x_i,x_j)\geq 2\sigma$ for $i\neq j$, and let $D_i$ be the Dirichlet-Voronoi cell of $x_i$. If
$\Xi$ is a quasi-orthoscheme associated to some $D_i$ and it is known that $\Xi$ is not an orthoscheme,  then
$$
\frac{|\Xi\cap B(x_i,\sigma)|}{|\Xi|\cdot |B(x_i,\sigma)|}
\leq \Delta(r_1(\sigma),\ldots,r_{d-2}(\sigma),r_{\infty}(\sigma)).
$$
\end{lemma}
\proof
Let $F_0\subset \ldots\subset F_{d-2}$ be the proper tower of faces of $D_i$ associated to $\Xi$.
If $\delta(x_i,q_i(F_{d-2}))\geq r_{\infty}(\sigma)$, then $F_{d-2}$ does not intersect the interior of
$B(x_i,r_\infty(\sigma))$, and hence Lemma~\ref{qiF} yields
$$
\frac{|\Xi\cap B(x_i,\sigma)|}{|\Xi|}\leq
\frac{|\Xi\cap B(x_i,\sigma)|}{|\Xi\cap B(x_i,r_\infty(\sigma))|}=
 \frac{|B(x_i,\sigma)|}{|B(x_i,r_\infty(\sigma))|}.
$$
Since $\Theta(r_1(\sigma),\ldots,r_{d-2}(\sigma),r_{\infty}(\sigma))\subset
B(z_0,r_\infty(\sigma))$, we have
$$
\frac{|\Theta(r_1(\sigma),\ldots,r_{d-2}(\sigma),r_{\infty}(\sigma))\cap B(z_0,\sigma)|}
{|\Theta(r_1(\sigma),\ldots,r_{d-2}(\sigma),r_{\infty}(\sigma))|}\geq
 \frac{|B(z_0,\sigma))}{|B(z_0,r_\infty(\sigma))|}.
$$
we conclude the lemma in this case.

This covers the case $d=3$ completely because the condition $\delta(x_i,q_i(F_1))< r_{\infty}(\sigma)$ implies by Lemma~\ref{qiF} that $\Xi$ is an orthoscheme. The only case left open is when $d=4$, $\delta(x_i,q_i(F_2))<r_{\infty}(\sigma)$, and hence $q_i(F_{2})$ is contained in the relative interior of $F_{2}$, but $q_i(F_{1})$ is not contained in the relative interior of $F_{1}$ because otherwise $\Xi$ is an orthoscheme. Then there exists $p\in\overline{q_i(F_{2}),q_i(F_1)}$ such that $\delta(x_i,p)=r_\infty(\varphi)$. We consider the spherical cone $C$ obtained by rotating the triangle with vertices $x_i,q_2(F_2),p$ around $\overline{x_i,q_2(F_2)}$. Since
$F_2\backslash C$ does not intersect $B(x_i,r_\infty(\varphi))$, the argument as above
leads to
\begin{equation}
\label{outside-cone}
\frac{|(\Xi\backslash C)\cap B(x_i,\sigma)|}
{|(\Xi\backslash C)|\cdot |B(x_i,\sigma)|}\leq
\Delta(r_1(\sigma),r_2(\sigma),r_{\infty}(\sigma)).
\end{equation}
In addition, (\ref{orthoscheme-density-monotone}) and the argument of K. Boroczky \cite{Bor78} yield
\begin{eqnarray}
\nonumber
\frac{|C\cap B(x_i,\sigma)|}
{|C|\cdot |B(x_i,\sigma)|}&=&\lim_{s\to 0^+}
\Delta(r_1(\sigma),r_{\infty}(\sigma)-s,r_{\infty}(\sigma))\\
\label{inside-cone}
&\leq&
\Delta(r_1(\sigma),r_2(\sigma),r_{\infty}(\sigma)).
\end{eqnarray}
Combining (\ref{outside-cone}) and (\ref{inside-cone}) proves Lemma~\ref{quasi-orthoscheme-density}.
\proofbox

Actually, the argument in K. Boroczky \cite{Bor78} shows that Lemma~\ref{quasi-orthoscheme-density} holds in any dimension. More precisely, \cite{Bor78} proved the so-called {\it simplex bound}; namely, if $\sigma\in(0,\frac{\pi}2)$,
and there exist $k$ non-overlapping spherical balls of radius $\sigma$ on $S^{d-1}$, then
\begin{equation}
\label{simplex-bound}
k\leq \Delta(r_1(\sigma),\ldots,r_{d-1}(\sigma))\cdot {\cal H}^{d-1}(S^{d-1}),
\end{equation}
and equality holds in the simplex bound if and only if the centers are vertices of a regular simplicial polytope $P$ with edge length $2\sin\sigma$.

The following statement shows in a qualitative way that if for an acute angle $\varphi$, all simplices in a Delone triangulation of $S^{d-1}$ are close to be regular with spherical edge length $2\varphi$, then the whole Delone triangulation is close to a one induced by a simplicial regular polytope.

\begin{lemma}
\label{close-simplices}
 Let $\varphi\in(0,\pi/4]$, let $u_0,\ldots,u_d\in S^{d-1}$, $d\geq 3$ be such that $u_1,\ldots,u_{d-1}$ determines a unique $(d-2)$-dimensional great subsphere that separates $u_0$ and $u_d$, and let $\varepsilon\in(0,\varepsilon_0)$ for $\varepsilon_0=\frac{\sin\varphi}{16\sqrt{d-1}}$.
If there exist two spherical regular simplices of edge length $\varphi$ with vertices
$v_0,\ldots,v_{d-1}$ and $w_1,\ldots,w_d$  such that $\delta(u_i,v_i)\leq \varepsilon$ for $i=0,\ldots, d-1$,
 and $\delta(u_i,w_i)\leq \varepsilon$ for $i=1,\ldots, d$,
then $\delta(u_d,v_d)\leq c\varepsilon$, where $v_1,\ldots,v_d$ are vertices of a regular simplex, $v_d\neq v_0$ and $c=\frac{16\sqrt{d-1}}{\sin\varphi}$.
\end{lemma}
\proof It is sufficient to prove that $\delta(v_d,w_d)\leq (c-1)\varepsilon$. Using $\delta(v_d,w_d) = 2\arcsin \frac {\|v_d-w_d\|} 2 \leq 2\|v_d-w_d\|$ given $\|v_d-w_d\|\leq 1$, it is sufficient to show
\begin{equation}
\label{close-simplices0}
\|v_d-w_d\|\leq \frac{c-1}2\cdot \varepsilon.
\end{equation}

We will use that if $x_1,\ldots,x_k,y_1,\ldots,y_k\in\R^k$, $\|x_i-y_i\|\leq\eta$ for all $i=1,\ldots,k$, and
$\lambda_1,,\ldots,\lambda_k\geq 0$, then the triangle inequality yields
\begin{equation}
\label{sum}
\|(\lambda_1x_1+\ldots+\lambda_kx_k)-(\lambda_1y_1+\ldots+\lambda_ky_k)\|
\leq (\lambda_1+\ldots+\lambda_k)\eta.
\end{equation}

We have $\delta(v_i,w_i)\leq 2\varepsilon$ for $i=1,\ldots, d-1$, thus
$\|v_i-w_i\|\leq 2\varepsilon$ for $i=1,\ldots, d-1$. We deduce from (\ref{sum}) that
$\|p-p'\|\leq 2\varepsilon$ holds for the centroids
$$
p=\frac1{d-1}(v_1+\ldots+v_{d-1})\mbox{ and }p'=\frac1{d-1}(w_1+\ldots+w_{d-1})
$$
of the $(d-2)$-dimensional regular Euclidean simplices $[v_1,\ldots,v_{d-1}]$ and $[w_1,\ldots,w_{d-1}]$.
We consider $\alpha>\beta>0$, and an orthonormal basis $\tilde{v}_1,\ldots,\tilde{v}_d$ such that
$v_d,\tilde{v}_d$ lie in the same half-space with respect to
${\rm lin}\{v_1,\ldots,v_{d-1}\}={\rm lin}\{\tilde{v}_1,\ldots,\tilde{v}_{d-1}\}$ and satisfy
\begin{eqnarray}
\label{vialphabeta}
v_i&=&\alpha\tilde{v}_i+\sum_{j\neq i\atop j\in\{1,\ldots,d-1\}}\beta\tilde{v}_j
\mbox{ \ for $i=1,\ldots,d-1$}
\end{eqnarray}
Then $\alpha,\beta$ satisfy
\begin{eqnarray*}
1&=&\langle v_1, v_1\rangle=\alpha^2+(d-2)\beta^2\\
\cos2\varphi&=&\langle v_1, v_2\rangle=2\alpha\beta +(d-3)\beta^2,
\end{eqnarray*}
therefore taking the difference leads to
\begin{equation}
\label{alphaminusbeta}
\frac{(\alpha-\beta)^2}2=\frac{1-\cos2\varphi}2=\sin^2\varphi.
\end{equation}
Similarly, we define an orthonormal basis $\tilde{w}_1,\ldots,\tilde{w}_d$ of $\R^d$  such that
$w_d,\tilde{w}_d$ lie in the same half-space with respect to
${\rm lin}\{w_1,\ldots,w_{d-1}\}={\rm lin}\{\tilde{w}_1,\ldots,\tilde{w}_{d-1}\}$ and satisfy
\begin{eqnarray}
\nonumber w_i&=&\alpha\tilde{w}_i+\sum_{j\neq i\atop j\in\{1,\ldots,d-1\}}\beta\tilde{w}_j
\mbox{ \ for $i=1,\ldots,d-1$}.
\end{eqnarray}
This basis exists when $\alpha,\beta$ satisfy the conditions derived above.

According to (\ref{vialphabeta}), the $(d-1)\times(d-1)$  symmetric matrix $M$ whose main diagonals are $\alpha$,
and the rest of the entries are $\beta$, satisfies that $M\tilde{v}_i=v_i$
$i=1,\ldots,d-1$. One of the eigenvectors of $M$ in $\tilde{v}^\bot_d$ is $v_*=\sum_{j=1}^{d-1}\tilde{v}_j$
with eigenvalue $\alpha+(d-2)\beta$. Any vector in $\tilde{v}^\bot_d$ orthogonal to $v_*$ is an eigenvector with eigenvalue $\alpha-\beta$. We deduce with help of (\ref{alphaminusbeta}) that if $v\in \tilde{v}^\bot_d$, then
\begin{equation}
\label{Minversenorm}
\|M^{-1}v\|\leq (\alpha-\beta)^{-1}\|v\|=\frac{\|v\|}{\sqrt{2}\sin\varphi}.
\end{equation}

For $i=1,\ldots,d-1$, we have $\langle\tilde{w}_d, w_i\rangle=0$ and $\|v_i-w_i\|\leq 2\varepsilon$, therefore,
$$
2\varepsilon\geq |\langle\tilde{w}_d, v_i\rangle|=\left\| \alpha\langle\tilde{w}_d,\tilde{v}_i\rangle+
\sum_{j\neq i\atop j\in\{1,\ldots,d-1\}}\beta\langle\tilde{w}_d,\tilde{v}_j\rangle \right\|.
$$
In particular, the length of the vector $v=\langle\tilde{w}_d, v_1\rangle \tilde{v}_1+\ldots+\langle\tilde{w}_d, v_{d-1}\rangle\tilde{v}_{d-1}$ is
at most $2\varepsilon\sqrt{d-1}$, thus (\ref{Minversenorm}) implies that
$$
\|M^{-1}v\| = \sqrt{\sum_{j=1}^{d-1}\langle\tilde{w}_d,\tilde{v}_j\rangle^2}\leq
\frac{2\varepsilon\sqrt{d-1}}{\sqrt{2}\sin\varphi}.
$$
In other words, the projection of the unit vector $\tilde{w}_d$ into $\tilde{v}^\bot_d$ is of length at most $\frac{2\varepsilon\sqrt{d-1}}{\sqrt{2}\sin\varphi}$, therefore, possibly after exchanging $\tilde{w}_d$ by $-\tilde{w}_d$, we have
$$
\|\tilde{v}_d-\tilde{w}_d\|\leq \frac{2\varepsilon\sqrt{d-1}}{\sqrt{2}\sin\varphi}\sqrt{2}=\frac{2\varepsilon\sqrt{d-1}}{\sin\varphi}.
$$
Now the orthogonal projection of the origin $o$ into ${\rm aff}\{v_1,\ldots,v_d\}$ lies inside $[p,v_d]$, thus the
angle of the triangle $[o,p,v_d]$ at $p$ is acute. In addition, the angle of $p$ and $v_d$ is also acute by
$\varphi\leq\frac{\pi}4$. Therefore, there exist $t,s\in(0,1)$ such that
$v_d=tp+s\tilde{v}_d$, and hence also $w_d=tp'+s\tilde{w}_d$. We deduce from $\|p-p'\|\leq 2\varepsilon\leq\frac{2\varepsilon\sqrt{d-1}}{\sin\varphi}$ and (\ref{sum}) that
$$
\|v_d-w_d\|\leq (t+s)\frac{2\varepsilon\sqrt{d-1}}{\sin\varphi}\leq \frac{4\varepsilon\sqrt{d-1}}{\sin\varphi}.
$$
According to (\ref{close-simplices0}), we may choose $c=\frac{16\sqrt{d-1}}{\sin\varphi}$.
\proofbox

We note that the lengthy calculations in the rest of paper (say, Section~\ref{secvol}) are mostly aiming at providing
upper estimates for the derivatives of $\Delta(\varphi_I-\varepsilon,r_2(\varphi_I-\varepsilon))$
 (see (\ref{icosa-Delta-epsilon})),
$\Delta(\varphi_I-\varepsilon,r_2(\varphi_I)+\gamma_2\varepsilon)$ (see Lemma~\ref{icosa-orhto-stab-long}),
$\Delta(\varphi_Q-\varepsilon,r_2(\varphi_Q-\varepsilon),r_3(\varphi_Q-\varepsilon))$ (see (\ref{600-Delta-epsilon}))
and $\Delta(\varphi_Q-\varepsilon,r_2(\varphi_Q-\varepsilon),r_3(\varphi_Q)+\gamma_3\varepsilon)$
(see Lema~\ref{600-orhto-stab-long}) as a function of small $\varepsilon>0$ where $\gamma_2$ and $\gamma_3$ are suitable large constants. These estimates can be obtained by some math computer packages based on formulas in
L. Fejes T\'oth \cite{FTL56} and \cite{FTL64}. However, we preferred a more theoretical approach, because the ideas can be used in any dimension for similar problems.

\section{Volume estimates related to the simplex bound}
\label{secvol}

To calculate or estimate $(d-1)$-volume of a compact $X\subset S^{d-1}$, we use Lemmas~\ref{radproj} and \ref{radtangent}.

\begin{lemma}
\label{radproj}
If $t\in(0,1)$, and $X\subset S^{d-1}$, $d\geq 3$, is spherically convex that for some $v\in X$ satisfies
$\langle u, v\rangle\geq t$ for all $u\in X$, then
${\cal H}^{d-1}(X)\geq {\cal H}^{d-1}(X')$ holds for the radial projection $X'$ of $X$ into $tv+v^\bot$.
\end{lemma}
\proof The statement follows from the fact that the orthogonal projection of $X$ into $tv+v^\bot$
covers $X'$. \proofbox

\begin{lemma}
\label{radtangent}
If $v\in S^{d-1}$, $d\geq 3$, and $X\subset S^{d-1}$ is compact and satisfies $\delta(u,v)\leq \Theta$, $\Theta<\frac \pi 2$, for all $u\in X$, and $\widetilde{X}$ is the radial projection of $X$ into the tangent hyperplane to $S^{d-1}$ at $v$, then
$$
{\cal H}^{d-1}(X)=\int_{\widetilde{X}}(1+\|y-v\|^2)^{-d/2}\,d{\cal H}^{d-1}(y)\geq \cos^d\Theta\cdot{\cal H}^{d-1}(\widetilde{X}).
$$
\end{lemma}
\proof The statement follows from the facts that if $y\in\widetilde{X}$, then
$\|y\|=(1+\|y-v\|^2)^{1/2}$ and $u=y/\|y\|$ satisfies $\langle u,v\rangle=(1+\|y-v\|^2)^{-1/2}\geq\cos\Theta$. \proofbox

The main results of this section are Lemma~\ref{orhto-volume-stab-regular}, its
Corollary~\ref{orhto-dense-stab-regular-coro}, and Lemma~\ref{orhto-volume-stab-long}, which provide estimates when
we slightly deform the "regular" orthoscheme $\Theta(r_1(\varphi),\ldots,r_{d-1}(\varphi))$.

\begin{lemma}
\label{orhto-volume-stab-regular}
For $\varphi\in \left(0,\arcsin \sqrt{\frac{d}{4(d-1)}}\right)$,  if $\varepsilon\in(0,\varphi)$, then
$$
|\Theta(r_1(\varphi-\varepsilon),\ldots,r_{d-1}(\varphi-\varepsilon))|>
|\Theta(r_1(\varphi),\ldots,r_{d-1}(\varphi))|(1-\aleph\cdot \varepsilon),
$$
where $\aleph=d\ 2^{(d+3)/2}/\sin r_{d-1}(\varphi)$.
\end{lemma}
\proof We deduce from (\ref{rjdefsin}) that  $r_{d-1}(\varphi)< \pi/4$.
Let $v\in S^{d-1}$, let $H=v+v^\bot$ be the hyperplane tangent to $S^{d-1}$ at $v$, and let $\sigma$ be a spherical arc of length $\pi/4$ starting from $v$. For $\varepsilon\in(0,\varphi)$, we consider the spherical regular simplex
$T(\varepsilon)$ whose spherical circumscribed ball is of center $v$ and radius $r_{d-1}(\varphi-\varepsilon)$, and one vertex of $T(\varepsilon)$ is contained in $\sigma$.
In particular,
$$
|\Theta(r_1(\varphi-\varepsilon),\ldots,r_{d-1}(\varphi-\varepsilon))|=|T(\varepsilon)|/d!.
$$
We write $\widetilde{T}(\varepsilon)$ to denote the radial projection of $T(\varepsilon)$ into $H$, which is a Euclidean regular simplex of circumradius
$R(\varepsilon)=\tan r_{d-1}(\varphi-\varepsilon)<1$.
Bounding ${\cal H}^{d-1}(\widetilde{T}(0))\leq 2^{\frac d2}|T(0)|$ by Lemma~\ref{radtangent} we deduce that
\begin{eqnarray}
\nonumber
|T(0)|-|T(\varepsilon)|&\leq& |\widetilde{T}(0)\backslash\widetilde{T}(\varepsilon)|\\
\nonumber
&=&\left(1-\frac{R(\varepsilon)^{d-1}}{R(0)^{d-1}}\right){\cal H}^{d-1}(\widetilde{T}(0))\\
\nonumber
&\leq&
\left(1-\left(1-\frac{R(0)-R(\varepsilon)}{R(0)}\right)^{d-1}\right)2^{d/2}|T(0)|\\
\label{R0Re}
&\leq&
\frac{R(0)-R(\varepsilon)}{R(0)}\cdot d\ 2^{d/2}|T(0)|.
\end{eqnarray}
For $r(\varepsilon)=r_{d-1}(\varphi-\varepsilon)$, we deduce from (\ref{rjdefsin}) that
$r'(\varepsilon)=-\frac{\cos(\varphi-\varepsilon)}{\cos r(\varepsilon)}\sqrt{\frac{2(d-1)}{d}}$, therefore,
$$
R'(\varepsilon)=(1+R(\varepsilon)^2)r'(\varepsilon)\geq -\frac{\sqrt{2}(1+R(\varepsilon)^2)}{\cos r(0)}
\geq -\frac{2^{3/2}}{\cos r(0)}.
$$
Using (\ref{R0Re}) and $R(0)\cdot \cos r(0)=\sin r_{d-1}(\varphi)$,
$$\frac {|T(0)|-|T(\varepsilon)|} {|T(0)|}\leq \frac {2^{3/2}\varepsilon}{R(0)\cdot \cos r(0)}\cdot d\ 2^{d/2}=\frac {d\ 2^{(d+3)/2}} {\sin r_{d-1}(\varphi)} \varepsilon.$$
\proofbox

\begin{coro}
\label{orhto-dense-stab-regular-coro}
For $\varphi\in \left(0,\arcsin \sqrt{\frac{d}{4(d-1)}}\right)$,  if $\varepsilon\in(0,\frac1{2\aleph})$ for the $\aleph$ of
Lemma~\ref{orhto-volume-stab-regular}, then
$$
\Delta(r_1(\varphi-\varepsilon),\ldots,r_{d-1}(\varphi-\varepsilon))\leq
\Delta(r_1(\varphi),\ldots,r_{d-1}(\varphi))(1+2\aleph\cdot \varepsilon).
$$
\end{coro}
\proof $1+2\aleph\varepsilon \geq 1/(1-\aleph\varepsilon)$ so, according to Lemma~\ref{orhto-volume-stab-regular}, it is sufficient to prove that if $0<s<\varphi$, then, for any $\tau<r_1(s)$,
\begin{equation}
\label{smallangle}
|B(z_0,\tau)\cap\Theta(r_1(s),\ldots,r_{d-1}(s))|
\leq |B(z_0,\tau)\cap\Theta(r_1(\varphi),\ldots,r_{d-1}(\varphi))|.
\end{equation}

Essentially, this statement means that the angle measure at a vertex of a regular spherical simplex increases when the side length of the simplex increases. For the sake of completeness we give an argument for this statement.

Consider two regular spherical simplices $T'$ and $T$ with side lengths $2s$ and $2\varphi$ respectively such that they share a common center $v$ and each vertex $z'_i$ of $T'$ belongs to the arc $\overline{z_i,v}$. Triangle $[z'_1,z'_2,v]$ is inside $[z_1,z_2,v]$ so the area of $[z'_1,z'_2,v]$ is less than the area of $[z_1,z_2,v]$. Since the area of a spherical triangle is the sum of its angles minus $\pi$, the angle between $\overline{z'_1,z'_2}$ and $\overline{z'_1,z'_v}$ is less than the angle between $\overline{z_1,z_2}$ and $\overline{z_1,v}$.

Now we consider two regular simplices $T'$ of side length $2s$ with vertices $z_0,z'_1,\ldots,z'_{d-1}$ and $T$ of side length $2\varphi$ with vertices $z_0,z_1,\ldots,z_{d-1}$ such that the center $v'$ of $T'$ belongs to the arc $\overline{v,z_0}$, where $v$ is the center of $T$, and all triangles $[z_0,v,z_i]$ and $[z_0,v',z'_i]$ overlap. Then all arcs $\overline{z_0,z_i}$ belong to the cone formed by $T$ at $z_0$ because all corresponding 2-dimensional angles in $T'$ are smaller than those in $T$. Therefore, the angle measure for $T'$ is smaller than the one for $T$.
\proofbox

We set up a notation for Lemma~\ref{orhto-volume-stab-long}. For $\varphi\in(0,\frac{\pi}4)$, let $z_0=z_0(\varphi)$, $z_1(\varphi)$, $\ldots$, $z_{d-1}(\varphi)$ be the vertices
of $\Theta(r_1(\varphi),\ldots,r_{d-1}(\varphi))$. For $t\in[r_{d-1}(\varphi),\frac{\pi}2)$, we
set
$$
\widetilde{\Theta}(\varphi,t)=\Theta(r_1(\varphi),\ldots,r_{d-2}(\varphi),t),
$$
and we
may assume that $z_0(\varphi),\ldots,z_{d-2}(\varphi)$ are vertices of $\widetilde{\Theta}(\varphi,t)$, and its $d$-th vertex $z_{d-1}(\varphi,t)$ satisfies
$z_{d-1}(\varphi)\in\overline{z_{d-2}(\varphi),z_{d-1}(\varphi,t)}$.

\begin{lemma}
\label{orhto-volume-stab-long}
If $\varphi\in \left(0,\arcsin \sqrt{\frac{d}{4(d-1)}}\right)$ and $t\in(\varphi,\frac{\pi}3)$, then
$$
\left|\widetilde{\Theta}(\varphi,t)\backslash \widetilde{\Theta}(\varphi,r_{d-1}(\varphi))\right|\geq
\frac{t-r_{d-1}(\varphi)}{2^d}
\cdot \left|\widetilde{\Theta}(\varphi,r_{d-1}(\varphi))\right|.
$$
\end{lemma}
\proof For brevity, we set $z_i=z_i(\varphi)$ for $i=0,\ldots,d-1$, and $r_{d-1}=r_{d-1}(\varphi)$. The condition on $\varphi$ yields that  $r_{d-1}\leq\frac{\pi}4$.

Let $s$ be the length of the arc $\overline{z_{d-1},z_{d-1}(\varphi,t)}$. Since the length of the
arc $\overline{z_{d-1},z_0}$ is $r_{d-1}$, and the angle of these two arcs is $\arccos\frac{-1}d$, the Law of Cosines (\ref{LawCosine}) yields
$$
\cos t=\cos r_{d-1} \cos s-(\sin r_{d-1} \sin s)/d,
$$
we deduce from $\sin t\geq \sin r_{d-1}$ that
$$
\frac{dt}{ds}=\frac{\cos r_{d-1} \sin s+(\sin r_{d-1} \cos s)/d}{\sin t}\leq \frac1{\sin r_{d-1}},
$$
therefore,
\begin{equation}
\label{st}
s\geq (t-r_{d-1})\sin r_{d-1}.
\end{equation}

We set $\widetilde{\Theta}=\widetilde{\Theta}(\varphi,r_{d-1}(\varphi))$, and observe that the closure of $\widetilde{\Theta}(\varphi,t)\backslash \widetilde{\Theta}$ is
the spherical simplex $T$ with vertices $z_0,\ldots,z_{d-3},z_{d-1},z_{d-1}(\varphi,t)$.
Let $H$ be the hyperplane tangent to $S^{d-1}$ at $z_{d-1}$, and  we write $X'$ to denote the radial projection of some $X\subset S^{d-1}$ in $H$. It follows that  $\widetilde{\Theta}'$ is the Euclidean orthoscheme such that
$d!$ of its copies tile the Euclidean regular simplex of circumradius $\tan r_{d-1}\leq 1$, and hence
$\|z'_{d-2}-z'_{d-1}\|=(\tan r_{d-1})/(d-1)$. We deduce from Lemma~\ref{radtangent} and (\ref{st}) that
\begin{eqnarray*}
|T|&\geq &\frac{|T'|}{2^d}=
\frac{|\widetilde{\Theta}'|\tan s}{2^d\|z'_{d-2}-z'_{d-1}\|}
\geq
\frac{|\widetilde{\Theta}'|(t-r_{d-1})\sin r_{d-1}}{2^d(\tan r_{d-1})/(d-1)}\\
&\geq& \frac{|\widetilde{\Theta}'|(t-r_{d-1})}{2^d}\geq
\frac{|\widetilde{\Theta}|(t-r_{d-1})}{2^d}.\mbox{ \ }\proofbox
\end{eqnarray*}

\section{The case of the Icosahedron}
\label{seciso}

In this section, we write $I$ to denote the regular icosahedron with vertices on $S^2$. In particular,
\begin{equation}
\varphi_I=\frac12\arccos\frac1{\sqrt{5}}<\arcsin \sqrt{\frac{3}{8}},
\end{equation}
thus Corollary~\ref{orhto-dense-stab-regular-coro} and Lemma~\ref{orhto-volume-stab-long} can be applied
with $\varphi=\varphi_I$.
Since $S^2$ can be dissected into $120$ congruent copies of $\Theta(\varphi_I,r_2(\varphi_I))$, we have
$$
|\Theta(\varphi_I,r_2(\varphi_I))|=\frac{\pi}{30},
$$
and it follows from (\ref{simplex-bound}) that
$$
\Delta(\varphi_I,r_2(\varphi_I))=\frac3{\pi}.
$$

According to (\ref{rjdefsin}), we have $\sin r_{2}(\varphi_I)=\frac2{\sqrt{3}}\sin\varphi_I$, thus
the constant $\aleph$ of
Lemma~\ref{orhto-volume-stab-regular} satisfies $\aleph=\frac{3\cdot 2^3}{\sin r_{2}(\varphi_I)}<40$.
In particular, Corollary~\ref{orhto-dense-stab-regular-coro} yields that if $\varepsilon\in(0,0.01)$, then
\begin{equation}
\label{icosa-Delta-epsilon}
\Delta(\varphi_I-\varepsilon,r_2(\varphi_I-\varepsilon))<\frac3{\pi}(1+80\varepsilon)
<\frac3{\pi}+80\varepsilon.
\end{equation}
We also note that if $v\in S^2$ and $\eta\in(0,\frac{\pi}2)$, then
\begin{equation}
\label{circle-area}
|B(v,\eta)|=2\pi(1-\cos\eta).
\end{equation}

\begin{lemma}
\label{icosa-orhto-stab-long}
For $\gamma\geq 10^4$ and $\varepsilon\in (0,\frac1{100\gamma})$, we have
$$
\Delta(\varphi_I-\varepsilon,r_2(\varphi_I)+\gamma\varepsilon)\leq
\Delta(\varphi_I,r_2(\varphi_I))-\frac{\gamma\varepsilon}{200}.
$$
\end{lemma}
\proof To simplify the notation, we write $\varphi=\varphi_I$ and
$r_2=r_2(\varphi)=\arcsin\frac{2\sin\varphi}{\sqrt{3}}$, which satisfy $r_2+\gamma\varepsilon<\frac{\pi}3$
(in order to apply Lemma~\ref{orhto-volume-stab-long}). We may assume that
$\Theta(\varphi-\varepsilon,r_2(\varphi-\varepsilon))$ and
$\Theta(\varphi-\varepsilon,r_2+\gamma\varepsilon)$ share a side of length $\varphi-\varepsilon$.

We deduce from $r_2(\varphi-\varepsilon)\leq r_2$ that $(r_2+\gamma\varepsilon) - r_2(\varphi-\varepsilon)\geq \gamma\varepsilon$.

%We deduce from (\ref{rjdefsin}) that $\frac{d}{d\varepsilon}\,r_2(\varphi-\varepsilon)=-\frac{2\cos(\varphi-\varepsilon)}{\sqrt{3}\cos r_2(\varphi-\varepsilon)}>-2$,and hence $r_2(\varphi-\varepsilon)\geq r_2-2\varepsilon$.

We set $T$ to be the closure of
$$
\Theta(\varphi-\varepsilon,r_2+\gamma\varepsilon)\backslash \Theta(\varphi-\varepsilon,r_2(\varphi-\varepsilon)),
$$
thus Lemma~\ref{orhto-volume-stab-long} yields
\begin{equation}
\label{dim2areaT}
|T|\geq
\frac{\gamma\varepsilon}{8}\cdot |\Theta(\varphi-\varepsilon,r_2(\varphi-\varepsilon))|.
\end{equation}
In addition, if $\sigma\in(0,\varphi-\varepsilon)$, then
we deduce from $\varepsilon<10^{-6}$, that
\begin{eqnarray*}
%\nonumber
\frac{|T\cap B(z_0,\sigma)|}{|B(z_0,\sigma)|\cdot |T|}&<&
\frac{|T\cap B(z_0,\sigma)|}{|B(z_0,\sigma)|\cdot |T\cap B(z_0,r_2(\varphi-\varepsilon))|}=
\frac{|B(z_0,\sigma)|}{|B(z_0,\sigma)|\cdot |B(z_0,r_2(\varphi-\varepsilon))|}\\
&\leq &
\frac1{|B(z_0,r_2(\varphi-10^{-6}))|}=\Delta_0<\frac3{\pi}-0.175,
\end{eqnarray*}

because $\Delta_0\approx 0.7751$ and $\frac3{\pi}-0.175\approx 0.7799$.

Therefore $\gamma\geq 10^4$ yields
\begin{eqnarray*}
\Delta(\varphi-\varepsilon,r_2+\gamma\varepsilon)&\leq&
\frac{(\frac3{\pi}+80\varepsilon)|\Theta(\varphi-\varepsilon,r_2(\varphi-\varepsilon))|+\Delta_0|T|}
{|\Theta(\varphi-\varepsilon,r_2(\varphi-\varepsilon))|+|T|}\\
&\leq &
\frac3{\pi}+80\varepsilon-\left(\frac3{\pi}+80\varepsilon-\Delta_0\right)\frac{\gamma\varepsilon/8}{1+\frac{\gamma\varepsilon}{8}}\\
&= &
\frac3{\pi}+ \gamma\varepsilon \left(\frac{80}{\gamma} -\frac {\frac3{\pi}+80\varepsilon-\Delta_0}{8+\gamma\varepsilon}\right)\\
&\leq &
\frac3{\pi}+ \gamma\varepsilon \left(10^{-2}-\frac {\frac3{\pi}-\Delta_0}{10}\right)\leq \frac3{\pi}-\frac{\gamma\varepsilon}{200}.
\end{eqnarray*}
\mbox{ }\hfill\proofbox

The following two simple statements are useful tools  in the case of the $600$-cell as well.

\begin{lemma}
\label{trianglearea}
If $T\subset \R^2$ is  a triangle such that all sides are of length at least $a$, and the center of the circle passing through the vertices lies in $T$, then $|T|\geq \frac{\sqrt{3}}4\,a^2$.
\end{lemma}
\proof The largest angle $\alpha$ of $T$ satisfies $\frac{\pi}3\leq\alpha\leq\frac{\pi}2$.
\proofbox

\begin{lemma}
\label{angle}
For $x,y,v\in S^2$, let $\delta(x,y)\geq2\psi$, and let $\delta(x,v)=\delta(y,v)=R$ for $0<\psi<R<\frac{\pi}2$.
If the angle between $\overline{v,x}$ and $\overline{v,y}$ is $\omega$, then
\begin{description}
\item{(i)} $\cos\omega\leq 1-\frac{2\sin^2\psi}{\sin^2R}$;
\item{(ii)}   If $\psi=\varphi-\varepsilon$ and $R\leq r+\gamma\varepsilon$
where $\psi<\varphi<r<\frac{\pi}2-\gamma\varepsilon$ and $\gamma>1$, then
$\cos\omega\leq 1-\frac{2\sin^2\varphi}{\sin^2 r}+\frac{4\gamma\varepsilon}{\sin^2r}$.
\end{description}
\end{lemma}
\proof For (i), the
Spherical Law of Cosines (\ref{LawCosine}) yields
$$
1-2\sin^2\psi=\cos 2\psi\geq
\cos^2 R+(\sin^2 R)\cos\omega=1-(1-\cos\omega)\sin^2R.
$$
Turning to (ii), we deduce from $\frac{d}{dt}\,\sin^2 t=\sin 2t\leq 1$ that
$$
\frac{2\sin^2(\varphi-\varepsilon)}{\sin^2 (r+\gamma\varepsilon)}\geq
\frac{2(\sin^2\varphi-\varepsilon)}{\sin^2 r+\gamma\varepsilon}=
\frac{(1-\frac{\varepsilon}{\sin^2\varphi})2\sin^2\varphi}{(1+\frac{\gamma\varepsilon}{\sin^2\varphi})\sin^2 r}
\geq \frac{\left(1-\frac{(\gamma+1)\varepsilon}{\sin^2\varphi}\right)2\sin^2\varphi}{\sin^2 r},
$$
and hence (i) implies (ii).
\proofbox

\noindent{\bf Proof of Theorem~\ref{simpl-pol-stability} in the case of the icosahedron }
Let $I$ be the icosahedron with vertices on $S^2$, therefore, the vertices determine the optimal packing of $12$ spherical circular discs of radius
$\varphi_I=\frac12\arccos\frac1{\sqrt{5}}$. We set $\varphi=\varphi_I$, $r_2=r_2(\varphi)$ and
$r_\infty=r_\infty(\varphi)$. For $\varepsilon_0=10^{-9}$ and $\eta=0.11$, we observe that
\begin{equation}
\label{r2eta}
r_2+10^7\varepsilon_0<r_2+\eta<r_\infty-\eta.
\end{equation}
Let $\varepsilon\in (0,\varepsilon_0)$, and let $x_1,\ldots,x_k\in S^2$ satisfy that $k\geq 12$, and
$\delta(x_i,x_j)\geq 2(\varphi-\varepsilon)$ for $i\neq j$. We may assume that for any $x\in S^2$ there exists $x_i$ such that
$\delta(x_i,x)<2(\varphi-\varepsilon)$. Let $P=[x_1,\ldots,x_k]$, and hence $o\in {\rm int}\, P$. We prove
Theorem~\ref{simpl-pol-stability} for the icosahedron in two steps.\\

\noindent {\bf Step 1 } {\it Proving that all Delone cells are of circumradius at most $r_2+10^7\varepsilon$}

We suppose that there exists a Delone cell of spherical circumradius at least
$r_2+10^7\varepsilon$, and seek a contradiction. Let us consider the triangulation of $S^2$ by all quasi-orthoschemes associated to the Dirichlet cell decomposition induced by  $x_1,\ldots,x_k$. Among them, let ${\cal O}$ and ${\cal Q}$ denote the family of the ones with diameter less than $r_2+10^7\varepsilon$, and with diameter at least $r_2+10^7\varepsilon$, respectively. We claim that
\begin{equation}
\label{dim2quasiarea}
\sum_{\Xi\in {\cal Q}}|\Xi|\geq 2\pi(1-\cos\eta)>0.03.
\end{equation}

Let $\varrho>0$ be the largest number such that $\varrho B^3\subset P$, and let $R=\arccos\varrho$. Then $\varrho B^3$ touches $\partial P$ at a point $y\in \partial P$ in the relative interior of a two-dimensional face $F$ of $P$, $R$ is the spherical circumradius of the corresponding Delone cell, and $R\geq r_2+10^7\varepsilon$. By construction, $R$ is the maximal circumradius among all Delone cells.

We may assume that $x_1,x_2,x_3$ are vertices of $F$ such that $y\in[x_1,x_2,x_3]=T$. Let $v=y/\|y\|$, and let $\widetilde{T}$ be the
radial projection of $T$ into $S^2$, that is the associated spherical "Delone triangle", and satisfies $v\in \widetilde{T}$.
If $R< r_\infty$, then all quasi-orthoschemes having vertex $v$ are actual orthoschemes by
Lemma~\ref{qiF}, and hence their union is $\widetilde{T}$. In particular, Lemmas~\ref{radproj} and \ref{trianglearea} yield that
$$
\sum_{\Xi\in {\cal Q}}|\Xi|\geq |\widetilde{T}|\geq |T|\geq
\frac{\sqrt{3}}4\,(2\sin(\varphi-\varepsilon_0))^2>0.4.
$$
However, if $R\geq r_\infty$ and $x\in B(v,\eta)$, then $\delta(x,x_i)\geq r_2+\eta$ for all $i=1,\ldots,k$,
thus any quasi-orthoscheme $\Xi$ containing $x$ has a diameter at least $r_2+10^7\varepsilon$ by (\ref{r2eta}). Therefore,
$$
\sum_{\Xi\in {\cal Q}}|\Xi|\geq |B(v,\eta)|=2\pi(1-\cos\eta)
$$
in this case, proving (\ref{dim2quasiarea}).

We note that $12=\frac3{\pi}\cdot|S^2|$ according to
the equality case of the simplex bound (\ref{simplex-bound}).
We deduce from (\ref{icosa-Delta-epsilon}), Lemma~\ref{icosa-orhto-stab-long} with $\gamma=10^7$
and (\ref{dim2quasiarea}) that
\begin{eqnarray*}
k&\leq &\sum_{\Xi\in {\cal O}}|\Xi|\frac3{\pi}\cdot(1+80\varepsilon)
+\sum_{\Xi\in {\cal Q}}|\Xi|\left(\frac3{\pi}-50,000\varepsilon\right)\\
&\leq& 12+\frac3{\pi}[4\pi\cdot 80\varepsilon- 0.03\cdot 50,000\cdot\varepsilon]<12.
\end{eqnarray*}
This contradiction completes the proof of Step 1.\\

\noindent{\bf Step 2 } {\it Assuming all Delone cells are of circumradius at most $r_2+10^7\varepsilon$}

It follows from (\ref{simplex-bound}) and (\ref{icosa-Delta-epsilon}) that $k=12$.

We set $\gamma=10^7$. Let $\Omega$ be a Delone cell, and let $v$ be the center of the circumcircle of radius $R$. We claim that $\Omega$ is a  triangle, and there exists a regular spherical triangle $\Omega_0$ of side length $2\varphi$, such that
for any vertex $x_i$ of $\Omega$ there exists a vertex $w$ of $\Omega_0$ with
\begin{equation}
\label{dim2Delone-regular}
\delta(x_i,w)\leq 25\gamma\varepsilon.
\end{equation}

If $x_i\neq x_j$ are the vertices of $\Omega$, and the angle between $\overline{v,x_i}$ and $\overline{v,x_j}$ is $\omega_{ij}$, then Lemma~\ref{angle},   $\sin\varphi/\sin r_2=\sqrt{3}/2$ and $\gamma\varepsilon<10^{-2}$ yield
$$
\cos\omega_{ij}\leq 1-\frac{2\sin^2\varphi}{\sin^2 r_2}+\frac{4\gamma\varepsilon}{\sin^2r_2}
\leq \frac{-1}2+12\gamma\varepsilon<0.
$$
In particular, $\Omega$ is a triangle by Corollary~\ref{d+2pack}.
Since $(\cos t)'=-\sin t$ is at most $\frac{-3}4$ if $t\in[\frac{\pi}2,\frac{2\pi}3]$, we have
\begin{equation}
\label{dim2gamma}
\omega_{ij}\geq\frac{2\pi}3-16\gamma\varepsilon.
\end{equation}
We deduce from the Remark after Theorem~\ref{simplex-unit-vectors} that one may find a regular spherical triangle $\Omega'$ with vertices on the spherical circle with center $v$ and radius $R$ such that
for any vertex $x_i$ of $\Omega$ there exists a vertex $w'$ of $\Omega'$ such that
the angle between $\overline{x_i,v}$ and $\overline{w',v}$ is at most $24\gamma\varepsilon$, and hence
$\delta(x_i,w')\leq 24\gamma\varepsilon$. We take $\Omega_0$ with the circumcenter $v$ so that for any vertex $w$ of $\Omega_0$ there exists a vertex $w'$ of $\Omega'$ such that $w\in\overline{w',v}$ or $w'\in\overline{w,v}$.
As $R\leq r_2+\gamma\varepsilon$ by the condition of Step 2, and $R\geq r_2(\varphi-\varepsilon)\geq r_2-\gamma\varepsilon$, we conclude (\ref{dim2Delone-regular})
by the triangle inequality.

 Now we fix a Delone cell $\Theta$ and let $\Theta_0$ be the spherical regular triangle provided by
(\ref{dim2Delone-regular}). We observe that $c<44$ for the constant of Lemma~\ref{close-simplices} in our case.
We may assume that the vertices of $\Theta_0$ are vertices of the face $F_0$ of the icosahedron $I$. There exist nine more faces $F_1,\ldots,F_9$ of $I$, such that  $F_i\cap F_{i-1}$ is a common edge for $i=1,\ldots,9$, and any vertex of $I$ is a vertex of some $F_i$, $i\leq 9$.
 Attaching the corresponding nine more Delone cells to $\Theta$, we conclude from Lemma~\ref{close-simplices} that we may choose
$c_I=44^9\cdot 25\gamma$.
\proofbox

\section{The case of the $600$-cell}
\label{sec600}

In this section, by $Q$ we denote the regular $600$-cell with vertices on $S^2$. In particular,
\begin{equation}
\varphi_Q=\frac{\pi}{10}<\arcsin \sqrt{\frac13}
\end{equation}
thus Corollary~\ref{orhto-dense-stab-regular-coro} and Lemma~\ref{orhto-volume-stab-long} can be applied
with $\varphi=\varphi_Q$.
Since $S^3$ can be dissected into $14400$ congruent copies of $\Theta(\varphi_Q,r_2(\varphi_Q),r_3(\varphi_Q))$, we have
$$
|\Theta(\varphi_Q,r_2(\varphi_Q),r_3(\varphi_Q))|=\frac{|S^3|}{14400}=\frac{\pi^2}{7200},
$$
and it follows from (\ref{simplex-bound}) that
\begin{equation}
\label{dim3Delta}
\Delta(\varphi_Q,r_2(\varphi_Q),r_3(\varphi_Q))=\frac{60}{\pi^2}.
\end{equation}
The main idea of the argument in the case of the $600$-cell will be similar to the one for the icosahedron.
According to (\ref{rjdefsin}), we have $\sin r_{3}(\varphi_Q)=\sqrt{\frac3{2}}\,\sin\varphi_Q$, thus
the constant $\aleph$ of
Lemma~\ref{orhto-volume-stab-regular} satisfies $\aleph=\frac{4\cdot 2^{7/2}}{\sin r_3(\varphi_Q)}<120$.
In particular, Corollary~\ref{orhto-dense-stab-regular-coro} yields that if $\varepsilon\in(0,0.004)$, then
\begin{equation}
\label{600-Delta-epsilon}
\Delta(\varphi_Q-\varepsilon,r_2(\varphi_Q-\varepsilon),r_3(\varphi_Q-\varepsilon))<\frac{60}{\pi^2}(1+240\varepsilon)<\frac{60}{\pi^2}+1500\varepsilon.
\end{equation}
Next Lemma~\ref{600-orhto-stab-long} estimates $\Delta(\varphi_Q-\varepsilon,r_2(\varphi_Q-\varepsilon),r_3(\varphi_Q)+\gamma\varepsilon)$ for large $\gamma$ and small $\varepsilon>0$, and Lemma~\ref{tetrahedron-volume} estimates the volume of a tetrahedron.

\begin{lemma}
\label{600-orhto-stab-long}
For $\gamma\geq 10^6$ and $\varepsilon\in (0,\frac1{100\gamma})$, we have
$$
\Delta(\varphi_Q-\varepsilon,r_2(\varphi_Q-\varepsilon),r_3(\varphi_Q)+\gamma\varepsilon)\leq
\Delta(\varphi_Q,r_2(\varphi_Q),r_3(\varphi_Q))-\frac{\gamma\varepsilon}{100}.
$$
\end{lemma}
\proof To simplify notation, we write $\varphi=\varphi_Q$ and
$r_3=r_3(\varphi)=\arcsin\frac{3\sin\varphi}{2}$,
and use the notation set up before Lemma~\ref{orhto-volume-stab-long}.

%We deduce from (\ref{rjdefsin}) that $\frac{d}{d\varepsilon}\,r_3(\varphi-\varepsilon)=-\frac{3\cos(\varphi-\varepsilon)}{2\cos r_3(\varphi-\varepsilon)}>-2$, and hence $r_3(\varphi-\varepsilon)\geq r_3-2\varepsilon$.

We deduce from $r_3(\varphi-\varepsilon)\leq r_3$ that $(r_3+\gamma\varepsilon) - r_3(\varphi-\varepsilon)\geq \gamma\varepsilon$.

For the closure $T$ of
$$
\widetilde{\Theta}(\varphi-\varepsilon,r_3+\gamma\varepsilon)\backslash
\widetilde{\Theta}(\varphi-\varepsilon,r_3(\varphi-\varepsilon)),
$$
 Lemma~\ref{orhto-volume-stab-long} yields
\begin{equation}
\label{dim3areaT}
|T|\geq
\frac{\gamma\varepsilon}{16}\cdot |\widetilde{\Theta}(\varphi-\varepsilon,r_3(\varphi-\varepsilon))|.
\end{equation}

Let $\sigma\in(0,\varphi-\varepsilon_0)$.
We consider two spherical cones $C$ and $C_0$, where
$C$ is obtained by rotating the triangle with vertices $z_0,z_1(\varphi-\varepsilon),z_3(\varphi-\varepsilon)$ around $\overline{z_0,z_1(\varphi-\varepsilon)}$, and
$C_0$ is obtained by rotating the triangle with vertices $z_0,z_1(\varphi-\varepsilon_0),
z_3(\varphi-\varepsilon_0)$ around $\overline{z_0,z_1(\varphi-\varepsilon_0)}$. For the two-face $F$ of $T$ opposite to $z_0$, $F\backslash C$ is disjoint from $B(z_0, r_3(\varphi-\varepsilon))$, which in turn contains $C$, and hence we have the density estimates
$$
\frac{|(T\backslash C)\cap B(z_0,\sigma)|}
{|T\backslash C|\cdot |B(z_0,\sigma)|}\leq
\frac{|B(z_0,\sigma)|}
{|B(z_0,r_3(\varphi-\varepsilon))|\cdot |B(z_0,\sigma)|}
\leq \frac{|C\cap B(z_0,\sigma)|}
{|C|\cdot |B(z_0,\sigma)|}.
$$
Since the density of $B(z_0,\sigma)$ in $C\cap T$ is
$\frac{|C\cap B(z_0,\sigma)|}{|C|}$, and in $T\backslash C$ the density
is at most $\frac{|C\cap B(z_0,\sigma)|}{|C|}$,
 we deduce using  (\ref{orthoscheme-density-monotone}) and the argument of K. Boroczky \cite{Bor78} that
\begin{eqnarray}
\nonumber
\frac{|T\cap B(z_0,\sigma)|}
{|T|\cdot |B(z_0,\sigma)|}&\leq&\frac{|C\cap B(z_0,\sigma)|}
{|C|\cdot |B(z_0,\sigma)|}
=\lim_{s\to 0^+}
\Delta(\varphi-\varepsilon,r_3(\varphi-\varepsilon)-s,r_3(\varphi-\varepsilon))\\
\nonumber
&\leq& \lim_{s\to 0^+}
\Delta(\varphi-\varepsilon_0,r_3(\varphi-\varepsilon_0)-s,r_3(\varphi-\varepsilon_0))\\
\label{smalldense}
&\leq &\frac{|C_0\cap B(z_0,\sigma)|}
{|C_0|\cdot |B(z_0,\sigma)|}=\Delta_0.
\end{eqnarray}

Now $C_0$ is a spherical cone whose base is a circular disc of radius
$\xi=\arccos\frac{\cos r_3(\varphi-\varepsilon_0)}{\cos (\varphi-\varepsilon_0)}$, center
$z_1(\varphi-\varepsilon_0)$ and height
$\varphi-\varepsilon_0$. Let $H\subset \R^4$ be the hyperplane tangent to $S^3$ at $z_1(\varphi-\varepsilon_0)$, let $C'_0$ be the radial projection of $C_0$ into $H$, which is a Euclidean cone whose base is a circular disc of radius
$\varrho=\tan \xi$, center  $z_1(\varphi-\varepsilon_0)$ and height
$h=\tan(\varphi-\varepsilon_0)$. Therefore, Lemma~\ref{radtangent} yields
\begin{eqnarray*}
|C_0|&=&\int_{C'_0}(1+\|x-z_1(\varphi-\varepsilon_0)\|^2)^{-2}\,dx\\
&=&
\int_0^{h}\int_0^{\varrho-\frac{\varrho t}h}(1+t^2+r^2)^{-2}\cdot 2\pi r\, drdt.
\end{eqnarray*}
In addition, if the angle between the arcs $\overline{z_0,z_1(\varphi-\varepsilon_0)}$ and $\overline{z_0,z_3(\varphi-\varepsilon_0)}$
is $\alpha$, then $\cos \alpha=
\frac{\tan (\varphi-\varepsilon_0)}{\tan r_3(\varphi-\varepsilon_0)}$. Therefore, (\ref{circle-area}) yields
$$
\Delta_0=\frac{1-\cos\alpha}{2|C_0|}<\frac{60}{\pi^2}-0.3.
$$
For $\Delta=\Delta(\varphi-\varepsilon,r_2(\varphi-\varepsilon),r_3+\gamma\varepsilon)$, $\gamma\geq 10^6$ yields
\begin{eqnarray*}
\Delta&\leq&
\frac{(\frac{60}{\pi^2}+1500\varepsilon)
\widetilde{\Theta}(\varphi-\varepsilon,r_3(\varphi-\varepsilon))|+\Delta_0|T|}{|\widetilde{\Theta}(\varphi-\varepsilon,r_2(\varphi-\varepsilon))|+|T|}\\
&\leq &
\frac{60}{\pi^2}+1500\varepsilon-
\left(\frac{60}{\pi^2}+1500\varepsilon-\Delta_0\right)\frac{\gamma\varepsilon/16}{1+\frac{\gamma\varepsilon}{16}}\\
&= &
\frac{60}{\pi^2}+\gamma\varepsilon \left(\frac {1500}{\gamma}-\frac
{\frac{60}{\pi^2}+1500\varepsilon-\Delta_0}{16+\gamma\varepsilon}\right)\\
&\leq &
\frac{60}{\pi^2}+\gamma\varepsilon \left(2\cdot 10^{-3}-\frac
{\frac{60}{\pi^2}-\Delta_0}{20}\right)\leq \frac{60}{\pi^2}-\frac{\gamma\varepsilon}{100}.
\end{eqnarray*}
\mbox{ }\hfill\proofbox

\begin{lemma}
\label{tetrahedron-volume}
If $\theta\in (0,\frac13)$, and $u_1,u_2,u_3,u_4\in S^2$ satisfy that $\langle u_i, u_j\rangle\leq -\theta$ for $i\neq j$, then
$$
{\cal H}^3([u_1,u_2,u_3,u_4])\geq \sqrt{\theta}/4.
$$
\end{lemma}
\proof For $T=[u_1,u_2,u_3,u_4]$, we have $o\in{\rm int}\ T$ by Lemma~\ref{inhemisphere}.
Let $r>0$ be the maximal number such that $rB^3\subset T$, and hence $r\leq \frac13$ (see, say, K. Boroczky \cite{Bor04}, Section~6.5).  We may assume that $rB^3$ touches $\partial T$ in a point $y$ of $F=[u_1,u_2,u_3]$, which lies in the relative interior of $F$. We set $u=y/r\in S^2$, and $v_i=(u_i-y)/\sqrt{1-r^2}\in S^2$ for $i=1,2,3$. We have
$\alpha\in[\arccos \frac13,\frac{\pi}2)$ and $\beta\in(\frac{\pi}2,\pi]$ such that
$\delta(u_i,u)=\alpha$ for $i=1,2,3$, $\delta(u_4,u)=\beta$. Thus $u_i=u\cos\alpha+v_i\sin\alpha$ for $i=1,2,3$, and
$u_4=-u|\cos\beta|+w\sin\beta $ for some $w\in u^\bot\cap S^2$.

Since $\langle u_i, u_j\rangle<0$ for $1\leq i<j\leq 3$, we have
$\langle v_i, v_j\rangle=\langle u_i, u_j\rangle-\cos\alpha\cos\alpha<0$ for $1\leq i<j\leq 3$.
We deduce that $\|u_i-u_j\|\geq \sqrt{2(1-r^2)}$ for $1\leq i<j\leq 3$, and
 there exists $l\in\{1,2,3\}$
such that $\langle v_l, w\rangle>0$. In particular, we have
$$
-\theta\geq \langle u_4, u_l\rangle\geq -|\cos \beta|\cdot \cos\alpha.
$$
It follows from Lemma~\ref{trianglearea} and $1-r^2\geq \frac89$ that
$$
{\cal H}^3(T)=\frac{|\cos \beta|+ \cos\alpha}4\cdot {\cal H}^2(F)\geq
\frac{\sqrt{|\cos \beta|\cdot\cos\alpha}}2\cdot \frac{\sqrt{3}(1-r^2)}2
>\frac{\sqrt{\theta}}4.
$$
\mbox{ }\hfill \proofbox

It is not hard to see that the lower bound $\sqrt{\theta}/4$ in Lemma~\ref{tetrahedron-volume} can't be replaced by, say, $2\sqrt{\theta}$. \\

\noindent{\bf Proof of Theorem~\ref{simpl-pol-stability} in the case of the $600$-cell }
Let $Q$ be an $600$-cell with vertices on $S^3$, therefore, its vertices determine the optimal packing of $120$ spherical circular discs of radius
$\varphi_Q=\frac{\pi}{10}$. We set $\varphi=\varphi_Q$, $r_2=r_2(\varphi)$, $r_3=r_3(\varphi)$ and
$r_\infty=r_\infty(\varphi)$. For $\gamma=10^{12}$, $\varepsilon_0=10^{-14}$ and $\eta=0.02$, we observe that
\begin{equation}
\label{r3eta}
r_3+\gamma\varepsilon_0<r_3+\eta<r_\infty-2\eta.
\end{equation}

Let $\varepsilon\in (0,\varepsilon_0)$, and let $x_1,\ldots,x_k\in S^2$ satisfy that $k\geq 120$, and
$\delta(x_i,x_j)\geq 2(\varphi-\varepsilon)$ for $i\neq j$. We may assume that for any $x\in S^3$,
there exists $x_i$ such that
$\delta(x_i,x)<2(\varphi-\varepsilon)$. Let $P=[x_1,\ldots,x_k]$, and hence $o\in {\rm int}\, P$. We prove
Theorem~\ref{simpl-pol-stability} for the $600$-cell in two steps.\\

\noindent {\bf Step 1 } {\it Proving that all Delone cells are of circumradius at most $r_3+\gamma\varepsilon$}

We suppose that there exists a Delone cell of spherical circumradius at least
$r_3+\gamma\varepsilon$ and seek a contradiction. Let us consider the triangulation of $S^3$ by all quasi-orthoschemes associated to the Dirichlet cell decomposition induced by  $x_1,\ldots,x_k$. Among them, let ${\cal O}$ and ${\cal Q}$ denote the family of the ones with diameter less than $r_3+\gamma\varepsilon$, and with diameter at least $r_3+\gamma\varepsilon$, respectively. We claim that
\begin{equation}
\label{dim3quasiarea}
\sum_{\Xi\in {\cal Q}}|\Xi|> (4\pi/3) \sin^3\eta>10^{-5}.
\end{equation}

Let $\varrho>0$ be the largest number such that $\varrho B^4\subset P$ and let $R=\arccos\varrho$. Then $\varrho B^4$ touches $\partial P$ at a point $y\in \partial P$ in the relative interior of a three-dimensional face $F$ of $P$, $R$ is the
spherical circumradius of the corresponding Delone cell, and
$R\geq r_3+\gamma\varepsilon$.

 We may assume that $x_1,x_2,x_3,x_4$ are vertices of $F$ in a way such that $y\in[x_1,x_2,x_3,x_4]=T$.
Let $v=y/\|y\|$, and let $\widetilde{T}$ be the
radial projection of $T$ into $S^3$, that is the associated spherical "Delone simplex", and satisfies $v\in \widetilde{T}$.
If $R< r_3+2\eta$, then all quasi-orthoschemes having vertex $v$ are actual orthoschemes by
Lemma~\ref{qiF}, and hence their union is $\widetilde{T}$.
If for some $\{i,j\}\subset\{1,2,3,4\}$, the angle between $\overline{v,x_i}$ and $\overline{v,x_j}$ is $\omega_{ij}$, then
Lemma~\ref{angle} yields
$$
\cos\omega_{ij}\leq 1-\frac{2\sin^2(\varphi-\varepsilon)}{\sin^2R}<
1-\frac{2\sin^2(\varphi-\varepsilon_0)}{\sin^2(r_3+2\eta)}<-0.1.
$$
In particular, Lemmas~\ref{radproj} and
\ref{tetrahedron-volume} yield that
$$
\sum_{\Xi\in {\cal Q}}|\Xi|\geq |\widetilde{T}|\geq |T|\geq \sqrt{0.1}/4>0.07.
$$
However, if $R\geq r_3+2\eta$ and $x\in B(v,\eta)$, then $\delta(x,x_i)\geq r_3+\eta$ for all $i=1,\ldots,k$,
thus any quasi-orthoscheme $\Xi$ containing $x$ has diameter at least $r_3+\gamma\varepsilon$ by (\ref{r3eta}).
We deduce from Lemma~\ref{radproj} that
$$
\sum_{\Xi\in {\cal Q}}|\Xi|\geq |B(v,\eta)|=(4\pi/3)\sin^3\eta
$$
in this case, proving (\ref{dim3quasiarea}).

We note that $120=\frac{60}{\pi^2}\cdot|S^3|$ according to
the equality case of the simplex bound (\ref{simplex-bound}).
We deduce from (\ref{icosa-Delta-epsilon}), Lemma~\ref{icosa-orhto-stab-long} with $\gamma=10^{12}$
and (\ref{dim2quasiarea}) that
\begin{eqnarray*}
k&\leq &\sum_{\Xi\in {\cal O}}|\Xi|\frac{60}{\pi^2}\cdot(1+1500\varepsilon)
+\sum_{\Xi\in {\cal Q}}|\Xi|\frac{60}{\pi^2}\cdot (1-10^{10}\cdot \varepsilon)\\
&\leq& 12+\frac{60}{\pi^2}[2\pi^2\cdot 1500\varepsilon- 10^{-5}\cdot 10^{10}\cdot\varepsilon]<12.
\end{eqnarray*}
This contradiction completes the proof of Step 1.\\

\noindent{\bf Step 2 } {\it Assuming all Delone cells are of circumradius at most $r_3+\gamma\varepsilon$}

It follows from (\ref{simplex-bound}) and (\ref{600-Delta-epsilon}) that $k=120$.

 Let $\Omega$ be a Delone cell, and let $v$ be the center of the circumscribed spherical ball of radius $R$. We claim that $\Omega$ is a  spherical tetrahedron and there exists a regular spherical tetrahedron $\Omega_0$ of side length $2\varphi$ such that
for any vertex $x_i$ of $\Omega$ there exists a vertex $w$ of $\Omega_0$ with
\begin{equation}
\label{dim3Delone-regular}
\delta(x_i,w)\leq 10,000 \gamma\varepsilon.
\end{equation}

If $x_i\neq x_j$ are the vertices of $\Omega$, and the angle between $\overline{v,x_i}$ and $\overline{v,x_j}$ is $\omega_{ij}$, then Lemma~\ref{angle},   $\sin\varphi/\sin r_3=\sqrt{2/3}$ and $\gamma\varepsilon<10^{-2}$ yield
$$
\cos\omega_{ij}\leq 1-\frac{2\sin^2\varphi}{\sin^2 r_2}+\frac{4\gamma\varepsilon}{\sin^2r_3}
\leq \frac{-1}3+30\gamma\varepsilon<0.
$$
In particular, $\Omega$ is a tetrahedron by Corollary~\ref{d+2pack}.
Since $(\cos t)'=-\sin t$ is at most $\frac{-3}4$ if $t\in[\frac{\pi}2,\frac{2\pi}3]$, we have
\begin{equation}
\label{dim3gamma}
\omega_{ij}\geq\arccos\frac{-1}3-40\gamma\varepsilon.
\end{equation}
We deduce from Theorem~\ref{simplex-unit-vectors} that one may find a regular spherical tetrahedron $\Omega'$
 with vertices on the subsphere with center $v$ and radius $R$ such that
for any vertex $x_i$ of $\Omega$ there exists a vertex $w'$ of $\Omega'$ such that
the angle between $\overline{x_i,v}$ and $\overline{w',v}$ is at most $9,000\gamma\varepsilon$ and hence
$\delta(x_i,w')\leq 9,000\gamma\varepsilon$. We take $\Omega_0$ with circumcenter $v$ so that for any vertex $w$ of $\Omega_0$ there exists a vertex $w'$ of $\Omega'$ such that $w\in\overline{w',v}$ or $w'\in\overline{w,v}$.
As $R\leq r_3+\gamma\varepsilon$ by the condition of Step 2, and $R\geq r_3(\varphi-\varepsilon)\geq r_3-\gamma\varepsilon$, we conclude (\ref{dim3Delone-regular})
by the triangle inequality.

 Now we fix a Delone cell $\Theta$ and let $\Theta_0$ be the spherical regular tetrahedron provided by
(\ref{dim3Delone-regular}). We observe that $c<90$ for the constant of Lemma~\ref{close-simplices} in our case.
We may assume that the vertices of $\Theta_0$ are vertices of the face $F_0$ of the $600$-cell $Q$. There exist $116$ more faces $F_1,\ldots,F_{116}$ of $Q$, such that  $F_i\cap F_{i-1}$ is a common edge for $i=1,\ldots,116$, and any vertex of $Q$ is a vertex of some $F_i$, $i\leq 116$.
 Attaching the corresponding $116$ more Delone cells to $\Theta$, we conclude from Lemma~\ref{close-simplices} that we may choose
$c_Q=90^{116}\cdot 10,000\,\gamma$.
\proofbox


\begin{thebibliography}{99}

\bibitem{Acz52}
J. Acz\'el:
Solution to Problem 35, I. (Hungarian).
Mat. Lapok, 3 (1952), 94-95.

\bibitem{And99}
N.N. Andreev:
A spherical code.
Russian Mathematical Surveys, 54 (1999), 251-253.

\bibitem{Bez10}
K. Bezdek:
Classical Topics in Discrete Geometry.
Springer, 2010.

\bibitem{BoF64}
K. B\"or\"oczky, A. Florian:
\"Uber die dichteste Kugelpackung im hyperbolischen Raum.
Acta Math. Hung., 15 (1964), 237-245.

\bibitem{Bor78}
K. B\"or\"oczky:
Packing of spheres in spaces of constant curvature.
Acta Math. Hung., 32 (1978), 243-261.

\bibitem{Bor04}
K. B\"or\"oczky, Jr.:
Finite packing and covering.
Cambridge Unversity Press, 2004.

\bibitem{BD01}
P. Boyvalenkov, D. Danev.
Uniqueness of the 120-point spherical 11-design in four dimensions.
Arch. Math., 77 (2001), 360-268.

\bibitem{CRY91}
S.A. Chepanov, S.S. Ryshkov, N.N. Yakovlev.
On the disjointness of point systems. (Russian)
Trudy Mat. Inst. Steklov., 196 (1991), 147-155.

\bibitem{CoS98}
J.H.~Conway, N.J.A.~Sloane:
Sphere packings,  Lattices and Groups.
Springer--Verlag, Berlin, New York 1998.

\bibitem{Cox73}
H.S.M. Coxeter:
Regular Polytopes. (originally published in 1947)
Dover, 1973.

\bibitem{DaH51}
H. Davenport, Gy. Haj\'os:
Problem 35 (Hungarian).
 Mat. Lapok, 2 (1951), 68.

\bibitem{ErZ01}
T. Ericson, V. Zinoviev:
Codes on Euclidean spheres.
North-Holland, 2001.

\bibitem{FTG04}
G. Fejes T\'oth:
Packing and covering.
In: E. J. Goodman, J. O´Rourke CRC Handbook on Discrete and Computational Geometry, CRC Press, 2004.

\bibitem{FTL49}
L. Fejes T\'oth:
On the densest packing of spherical caps.
Amer. Math. Monthly, 56 (1949), 330-331.

\bibitem{FTL56}
L. Fejes T\'oth:
On the volume of a polyhedron in non-Euclidean spaces.
Publ. Math. Debrecen 4 (1956), 256-261.

\bibitem{FTL64}
L. Fejes T\'oth:
Regular Figures.
Pergamon Press,
Oxford, 1964.



\bibitem{FTL72}
L. Fejes T\'oth:
Lagerungen in der Ebene, auf der Kugel und im Raum.
Springer, Berlin, 1972.

\bibitem{Jun01}
H.W.E. Jung:
\"Uber die kleineste Kugel die eine r\"aumliche Figur einschliesst.
J. Reine ang. Math., 123 (1901), 241-257.

\bibitem{McS02}
P. McMullen, E. Schulte:
Abstract Regular Polytopes.
Cambridge University Press. 2002.

\bibitem{Mus08}
O. Musin:
The kissing number in four dimensions.
Annals of Mathematics, 168 (2008), 1-32.

\bibitem{Ran55}
R.A. Rankin:
The closest packing of spherical caps in $n$ dimensions.
 Proc. Glasgow Math. Assoc., 2 (1955), 139-144.

\bibitem{Rog64}
 C.A. Rogers:
Packing and Covering.
Cambridge Univ.~Press, Cambridge, 1964.

\bibitem{Sze52}
T. Szele:
Solution to Problem 35, II. (Hungarian).
Mat. Lapok, 3 (1952), 95.

\end{thebibliography}
\end{document}